\documentclass[10pt]{article}

\usepackage[letterpaper, margin=0.75in]{geometry}
\usepackage{amsmath}
\usepackage{amssymb}
\usepackage{amsthm}
\usepackage{nomencl}
\usepackage{graphicx}
\usepackage{caption}
\usepackage{subcaption}
\usepackage{wrapfig}
\usepackage{multirow}
\usepackage{placeins}
\usepackage{hyperref}
\usepackage{algorithm}
\usepackage{algpseudocode}
\usepackage[backend=biber,bibencoding=utf-8,style=numeric,sorting=none]{biblatex}

\let\bld\boldsymbol

\graphicspath{{../../images/}{../../images/results-naca0012-abilu/}{images/}}

\newtheorem{theorem}{Theorem}
\newtheorem{lemma}{Lemma}
\newtheorem{definition}{Definition}

\title{An asynchronous incomplete block LU preconditioner for computational fluid dynamics on unstructured grids}
\author{Aditya Kashi and Siva Nadarajah}

\begin{document}

\maketitle

\begin{abstract}
We present a study of the effectiveness of asynchronous incomplete LU factorization preconditioners for the time-implicit solution of compressible flow problems while exploiting thread-parallelism within a compute node. A block variant of the asynchronous fine-grain parallel preconditioner adapted to a finite volume discretization of the compressible Navier-Stokes equations on unstructured grids is presented, and convergence theory is extended to the new variant. Experimental (numerical) results on the performance of these preconditioners on inviscid and viscous laminar two-dimensional steady-state test cases are reported. It is found, for these compressible flow problems, that the block variant performs much better in terms of convergence, parallel scalability and reliability than the original scalar asynchronous ILU preconditioner. For viscous flow, it is found that the ordering of unknowns may determine the success or failure of asynchronous block-ILU preconditioning, and an ordering of grid cells suitable for solving viscous problems is presented.

\end{abstract}


\section{Introduction}
Today's new high-performance computers have multiple levels of parallelism, from distributed memory parallelism of clusters to vector operations in individual processor cores. To further speed up computational fluid dynamics (CFD) codes and enable solutions with greater detail, these levels of parallelism need to be fully harnessed. In CFD, parallelism at the level of a cluster (a networked set of independent processors called nodes) is fairly mature. However, parallelism within each node of a cluster is a work-in-progress, especially for implicit methods of integration in time. With the development and commercialization of fine-grain parallel processors such as graphics processing units (GPUs) \cite{gpu:nvidiatesla} and many-core central processing units (CPUs) \cite{mic:knl}, node-level parallelism has become increasingly important.

In case of explicit time-stepping, there has been some success in utilizing many-core devices and GPUs (eg. \cite{ho:pyfr}). However, explicit solvers face a restrictive time-step limit and are only suitable when resolution of unsteady high-frequency phenomena is desired. 
For steady-state problems, implicit time-stepping is more suitable.
However, implicit solvers require the solution of large sparse systems of linear equations, which in turn, requires effective preconditioners.

Incomplete LU factorization is commonly used as a preconditioner to solve large sparse systems of equations due to its wide applicability. Several variants are available depending on parameters such as level-of-fill, threshold for a drop-tolerance and ordering of the grid points \cite{saad}.
Level of fill refers to allowing more non-zeros in the lower and upper triangular factors than in the original matrix at locations determined exclusively by the sparsity pattern of the original matrix. On the other hand, sparsity of the factors can be controlled depending on the magnitude of the introduced non-zeros and a drop tolerance. These two approaches can be combined. In addition, since the effectiveness of an ILU preconditioner depends on the ordering of the grid (and thus of the original matrix), the ordering can be changed to attempt to improve convergence.
However, factorization of the matrix into upper and unit lower triangular factors, and the application of the triangular factors during the solve, are originally inherently sequential. When only coarse-grained parallelism is necessary, sequential ILU can be used as a subdomain iteration for a domain decomposition preconditioner \cite{dd:gropp}. However, when fine grain parallelism is required, one first turns to the available parallel variants of ILU factorization and triangular solves.

A typical way of formulating a parallel ILU algorithm (and parallel triangular solves) is to re-order the unknowns and equations such that in the new ordering, some rows or columns can be eliminated in parallel.
One such ordering is the multi-colour ordering \cite{ilu:orderings:multicolour_iccg, ilu:orderings:parallel_benzi}. This involves labelling every unknown with a colour such that unknowns labelled with the same colour do not have a direct dependence on each other. The unknowns are then re-ordered such that those assigned the same colour are numbered consecutively, which leads to a block structure with (usually large) diagonal blocks.
The advantage of this approach is that the degree of parallelism is proportional to the problem size for common discretizations. However, researchers have observed \cite{ilu:orderings:multicolour_iccg, ilu:orderings:parallel_benzi} significant increases in the number of linear solver iterations compared to other non-parallel orderings. Furthermore, a multi-colour triangular solver may lead to bad memory access patterns, especially on GPUs, due to a non-unit stride \cite{multicolour:nguyen_gpu_2019}. This could hurt performance.

Another approach is level-scheduling. In this approach, an analysis phase computes (without necessarily re-ordering the unknowns) sets of rows that can be eliminated in parallel due to the sparsity of the matrix \cite{gs:naumov, ilu:naumov}. On structured grids with cells ordered by $i,j,k$, this approach gives rise to wavefront algorithms, where planes of cells having the same value of $i+j+k$ are processed in parallel \cite{ilu:edwards}. This method has the advantage of not adversely affecting the convergence rate while enabling reasonably good memory access patterns. However, the parallelism depends on the sparsity pattern of the matrix, and may be very limited in some cases such as the hybrid unstructured grids used in this work. In the case of a structured square grid in two dimensions with a total of $n$ cells, at most $\sqrt{n}$ cells can be processed in parallel at a time.
Level-scheduling and multi-colouring have been compared, eg. by Suchoski et al. \cite{gs:raghavan} in case of symmetric positive-definite matrices. In that work it was found that multi-colouring was generally faster than level-scheduling on GPUs, though the speedup was highly dependent on the matrix.
Another notable technique is a variant of level-scheduling designed for triangular systems \cite{gs:optimal_alvardo}, in which the matrix is reordered to compute an optimal partition of triangular factors. Each of the triangular factors can be inverted in-place in parallel. For several kinds of triangular matrices, this method exposes more parallelism than a typical level-scheduled method, though the amount of parallelism is seen to slowly decrease with problem size.


The idea of asynchronous iterations has been proposed as a way of devising parallel iterations. Originally referred to as chaotic relaxation by Chazan and Miranker \cite{async:chazan_1969}, these schemes aim to remove all synchronization in linear iterations and still obtain convergence. The concept has been extended to non-linear iterations as well \cite{async:frommer_2000}. An asynchronous iterative method to compute incomplete LU (ILU) factorization has been proposed by Chow and Patel \cite{ilu:chowpatel} and applied to discretized linear partial differential equations (PDEs). This technique is the basis for the extensions and applications presented here.

There have been some efforts to apply asynchronous iterations to solution of scalar partial differential equations (PDEs). In the 1980's, Anwar and El Tarazi used a non-linear asynchronous iteration to solve a Poisson problem with non-linear boundary conditions \cite{async:nonlinear_poisson}. More recently, Chow and Patel \cite{ilu:chowpatel} demonstrated their asynchronous ILU factorization for solving the Poisson equation and the linear convection-diffusion equation with promising results.
That being said, application of asynchronous iterations specifically to fluid dynamics problems has been rare. Chaotic relaxation was recently applied to incompressible flows in marine engineering \cite{async:hawkes_incompressible}. The solver used a pressure-based method whereby the momentum and mass flux equations are solved in separate steps. The authors demonstrate better strong scaling for chaotic relaxation compared to a Jacobi-preconditioned generalized minimum residual (GMRES) solver for certain problem sizes and a specific solution strategy. However, for the tightly coupled system of PDEs of compressible viscous flow, we show that regular chaotic relaxation and asynchronous ILU factorization may be insufficient. Second, the current authors and collaborators demonstrated the applicability of `asynchronous block symmetric Gauss-Seidel' (SGS) iterations as a fine-grain parallel multigrid smoother \cite{async:aditya_absgs_2019}. That work was in the context of structured grids only, but showed that good parallel scalability could be obtained on multi-core CPUs as well as the Intel Xeon Phi Knights Landing many-core processor, for viscous compressible flow problems in external aerodynamics.

This paper is concerned with fine-grain parallel subdomain preconditioners for implicit solvers of steady-state compressible flow problems on unstructured grids. By `subdomain preconditioner', we refer to the iteration applied locally on each subdomain of a domain decomposition preconditioner; the latter is used to parallelize across the nodes of a cluster. For this work, the `global' domain decomposition preconditioner (such as additive Schwartz) is not the subject of study. We propose a point-block variant of Chow and Patel's asynchronous ILU(0) factorization method \cite{ilu:chowpatel} as the subdomain preconditioner.

We also propose a point-block asynchronous iteration for applying the $L$ and $U$ factors.
In this context, Chow et al. \cite{triangular:chow_async_block_2018} showed that while using Jacobi iterations to apply triangular factors, it is advantageous to find blocks in the matrix and invert the diagonal blocks exactly. However, this was done only in the context of symmetric positive-definite matrices. In the proposed asynchronous approach, the applied preconditioner may change from one linear iteration to another. We therefore use flexible GMRES \cite{gmresflexible} as our Krylov subspace solver.
Additionally, we describe orderings of grid cells for effective solution of viscous problems by asynchronous block ILU(0) preconditioning. Finally, we show some experimental results for both inviscid and viscous compressible flow problems (section \ref{sec:results}), highlighting the necessity of the block variant and the effect of re-ordering on the performance of asynchronous ILU preconditioners.

\section{The solver setup}
\label{sec:solver}

In this work, we solve the compressible Navier-Stokes equations expressed in terms of the conserved variables density, momentum per unit volume in each spatial direction and total energy per unit volume as described in equations 2.19 through 2.24 by Blazek \cite[chapter~2]{blazek}.
To close the system and obtain $d+2$ equations for the $d+2$ conserved variables (in the case of $d$ spatial dimensions), the ideal gas equation of state, Fourier's law of heat conduction and Sutherland's law relating viscosity to temperature are used.

The steady-state equations are discretized in space to obtain a nonlinear system of equations
$\bld{r}(\bld{w}) = \bld{0}$,
where $\bld{w}\in \mathbb{R}^n$ with $n$ being the number of cells times the number of conserved variables ($d+2$) and $\bld{r}$ is the vector of the discretized residual functions on all cells stacked in the same order as $\bld{w}$.
This could be solved by the Newton-Raphson method or using `pseudo-time' stepping. In the latter approach, the discretized equations are converted into a nonlinear system of ordinary differential equations 
$\bld{M}\frac{d\bld{w}}{d\tau} + \bld{r}(\bld{w}) = \bld{0}$,
where $\bld{M}$ is a mass matrix depending on the spatial discretization. This is the approach implemented in our in-house CFD solver and used in this work. We can now iterate to steady-state after starting from a trivial initial condition such as uniform flow. We use a backward Euler discretization for the pseudo-time term:
\begin{align}
\left(\bld{M}_\tau^n + \frac{\partial\bld{r}}{\partial\bld{w}}(\bld{w}^n)\right)\Delta\bld{w}^n = -\bld{r}(\bld{w}^n), \label{eq:pseudotimestep} \\
\bld{w}^{n+1} = \bld{w}^n + \omega\Delta\bld{w}^n.
\end{align}
Above, $\omega \in [0.2,1)$ is a relaxation factor meant to prevent very large relative changes of density or pressure, and $\bld{M}_\tau^n$ is a diagonal matrix with diagonal blocks given by
$(\bld{M}_\tau^n)_i := \frac{V_i}{\Delta\tau_i^n}\bld{I}_{d+2}$,
where $V_i$ is the volume of the $i$th cell. The pseudo-time steps $\Delta\tau_i^n$ are determined using a CFL (Courant-Friedrichs-Lewy) number. The CFL number starts at a relatively small prescribed value and is ramped exponentially with respect to the ratio of residual norms from one time step to the next:
\begin{equation}
c_{fl}^{n+1} := c_{fl}^n \left( \frac{\lVert \bld{r}^n \rVert}{\lVert \bld{r}^{n+1} \rVert} \right)^r, \quad \text{where} \,\,\,
r := \begin{cases}
0.25 \quad \text{if}\,  \frac{\lVert \bld{r}^n \rVert}{\lVert \bld{r}^{n+1} \rVert} > 1 \\
0.3 \quad \text{if}\,  \frac{\lVert \bld{r}^n \rVert}{\lVert \bld{r}^{n+1} \rVert} \leq 1
\end{cases}.
\label{eqn:cfl}
\end{equation}
The exponents above were chosen empirically as they worked well for a number of test cases. The CFL number is kept between a prescribed starting minimum value and a maximum value.
At every pseudo-time step, the same CFL number is used for every cell. The local time step at cell $K_i$ with index $i$ and pseudo-time step $n$ is computed as
\begin{equation}
\Delta\tau^n_i = c_{fl}^n \,\,V_i\,/\sum_{f_{ij} \subset \partial K_i} \nu_{ij}(\bld{w}^n)\, a_{ij},
\end{equation}
where the sum is over the faces $f_{ij}$ of cell $K_i$, $\nu_{ij}$ is the absolute maximum eigenvalue of the normal flux Jacobian (including viscous fluxes) at face $f_{ij}$ and $a_{ij}$ is the area of that face.

The spatial discretization is a cell-centred finite volume scheme over two-dimensional unstructured hybrid grids. Upwind numerical fluxes are used for the inviscid terms. Gradients at cell-centres are estimated by a least-squares approach using data from face-neighbouring cells and reconstruction to faces is done by linear interpolation. Limiters were not used for the results shown here. Gradients at faces are computed as a modified average of the left and right cell-centred gradients. The Jacobian matrix used in \eqref{eq:pseudotimestep} is computed ignoring the reconstruction; that is, the first-order inviscid numerical flux and `thin-layer' first-order viscous flux are linearized to compute the Jacobian.

We note here that in our implementation, the ordering of the conserved variables $\bld{w}$ and of residuals in $\bld{r}$ are such that all $d+2$ variables of one cell are placed consecutively, followed by those of the next cell, and so on. For example in two dimensions,
\begin{equation}
\bld{w}^T = [\rho_1, \rho v_{x1}, \rho v_{y1}, \rho E_1, \,\, \rho_2,  \rho v_{x2}, \rho v_{y2}, \rho E_2,\, ...] = [\bld{u}_1^T \, \bld{u}_2^T \, ... \, \bld{u}_N^T],
\end{equation}
where subscripts denote cell indices. The residuals are ordered in the same way, with the mass, momentum and energy fluxes for a particular cell placed consecutively. This ordering leads to a block structure of the Jacobian with small dense blocks.
Preconditioners for such a blocked matrix are referred to as point-block preconditioners. Point blocking is advantageous for systems of PDEs because the small dense blocks can be inverted exactly to resolve the local coupling between the different physical variables at one mesh location. Once the required operations on blocks (inversion, matrix-vector products etc.) are available, this allows an extension of iterations for scalar PDEs to effectively deal with systems.
This is our main motivation for deriving block variants of asynchronous preconditioners.

The question of storage layout is orthogonal to the above discussion on point-block matrices. We store the matrix such that the entries in a block are stored contiguously. In case of point-block solvers for unstructured grids on many-core CPU architectures, this is advantageous for cache-locality and could be advantageous for vectorization (depending on the grid connectivity and block size). The proposed asynchronous block ILU preconditioner can also be used with different layouts suitable for other devices such as GPUs, but this is not discussed in this paper.


\nomenclature[01]{$v_j$}{The $j$th entry of vector $\bld{v}$} %
\nomenclature[02]{$\bld{v}_j$}{The $j$th sub-vector (of some specified size) of vector $\bld{v}$} %
\nomenclature[03]{$A_{ij}$}{The $(i,j)$th (scalar) entry of matrix $\bld{A}$} %
\nomenclature[04]{$\bld{A}_{ij}$}{The block at the $(i,j)$th block-index of matrix $\bld{A}$ according to some specified blocking} %

We consider the problem of parallel preconditioning of the system of equations given in \eqref{eq:pseudotimestep}, which we rewrite for simplicity as $\bld{A}\bld{x} = \bld{b}$,
where $\bld{A} \in \mathbb{R}^{n\times n}, \bld{x},\bld{b} \in \mathbb{R}^n$. For this work, we aim to replace this with the equivalent system $\bld{A}\bld{M}^{-1}\bld{M}\bld{x} = \bld{b}$
and choose the right-preconditioning matrix $\bld{M}$ such that this new system is better conditioned. While left preconditioning can also be used, we focus on right-preconditioning because the FGMRES solver requires it. As mentioned in the introduction, flexible Krylov subspace solvers are needed to accommodate variable preconditioners.

Incomplete LU (ILU) factorizations are commonly used as (sequential) preconditioners. ILU preconditioners are of the form $\bld{M}_{ilu} = \bld{L}\bld{U}$,
where $\bld{L}$, a unit lower triangular matrix and $\bld{U}$, an upper triangular matrix, approximate the LU factorization of $\bld{A}$.
Block ILU factorization is common in the solution of compressible flows by finite volume methods. Block factorization is defined here with respect to dense square blocks of size $b \times b$, such that $n$ is divisible by $b$.
While the results presented here can be extended to more general and variable block sizes, we adopt fixed-size square blocks for this work. For our problems, $b = d+2$ and each block represents the coupling between the $d+2$ fluxes and flow variables at two grid locations. Our finite volume discretization of the compressible Navier-Stokes equations gives dense, disjoint blocks. 
The matrix $\bld{A}$ and the block-triangular matrices $\bld{L}$ and $\bld{U}$ can now be viewed as being made up of these blocks and each block can be given an $(i,j)$ block-index depending on its location in the matrix with respect to other blocks. A $b\times b$ block is referred to as a `non-zero block' if there is at least one non-zero entry in the block.
For the block ILU factorization, $\bld{L}$ is a lower block triangular matrix with identity matrices for the diagonal blocks and $\bld{U}$ is an upper block triangular matrix with non-singular diagonal blocks.

Traditional algorithms for computing ILU preconditioners (eg., \cite[algorithm~10.1]{saad}) and applying them are sequential and difficult to parallelize; there is data dependency between tasks (eg., algorithm \ref{alg:ftrisubs}). As explained in the introduction, known parallel ILU and triangular solution schemes such as multi-colouring and level scheduling have drawbacks.
A parallel alternative can be derived using the concept of chaotic or asynchronous relaxation. We discuss this in the next section.
\begin{algorithm}[h]
	\caption{Forward substitution $\bld{L}\bld{x} = \bld{b}$} 
	\label{alg:ftrisubs}
	\begin{algorithmic}[1]
		\Require $\bld{L}$ is lower triangular, $\bld{b}$ is the right-hand side vector
		\State $x_1 \gets b_1/L_{11}$
		\For {$i$ from $2$ to $n$}
		\State $x_i := (b_i - \sum_{j=1}^{i-1} L_{ij}x_j)/L_{ii}$ \Comment Data dependency!
		\EndFor
	\end{algorithmic}
\end{algorithm}

\section{Review of asynchronous iterations}

\subsection{Chaotic relaxation}
\label{sec:chaotic}
Chazan and Miranker in 1969 had suggested \cite{async:chazan_1969} an asynchronous linear fixed-point iteration to solve the linear system $\bld{Ax}=\bld{b}$. They called their method `chaotic relaxation'; it is an asynchronous generalization of the well-known Jacobi and Gauss-Seidel iterations.

Let $\bld{A}$ be an $n\times n $ non-singular matrix and consider a splitting $\bld{A} = \bld{M}+\bld{N}$ defining the relaxation
\begin{equation}
\bld{x}^{j+1} = \bld{M}^{-1}\bld{b} - \bld{M}^{-1}\bld{N}\bld{x}^j.
\label{stationary_iteration}
\end{equation}
Let $\bld{B} := -\bld{M}^{-1}\bld{N}$ and $\bld{C} := \bld{M}^{-1}$ with rows $\bld{c}_i^T$. Chazan and Miranker define a chaotic relaxation as the following fixed-point iteration scheme.
\begin{equation}
x_i^{j+1} = 
\begin{cases} 
x_i^j & \text{if } i \neq u(j) \\
\sum_{\alpha=1}^{n} B_{i\alpha}x_\alpha^{j-s_\alpha(j)} + \bld{c}_i^T\bld{b} & \text{if } i = u(j).
\end{cases}
\label{eq:chaotic_relaxation}
\end{equation}
Following Strikwerda \cite{async:strikwerda_convergence_1997}, we refer to $s_\alpha:\mathbb{N}\rightarrow\mathbb{N},\, \alpha \in \{1,2,...,n\}$ as the `shift' or `delay' functions and $u:\mathbb{N}\rightarrow\{1,2,...,n\}$ as the `update function'.
Note that $j$ here does not refer to the $j$th iterate as in \eqref{stationary_iteration}, but a `step'. Frommer and Szyld consider a `step' be determined by one read of the approximate solution vector $\bld{x}$ from memory by one processing element (core, vector unit lane, CPU etc.) \cite{async:frommer_2000}. This is the convention we adopt here. The following conditions are imposed:
\begin{itemize}
	\item The shifts are bounded above uniformly:
	\begin{equation} 
	    \exists\, \hat{s} \in \mathbb{N}\,\, \text{s.t.}\,\, 0 \leq s_i(j) \leq \min \{j-1,\hat{s}\} \quad \forall i \in \{1,2,...n\},\, j \in \mathbb{N}.
	    \label{eq:chaotic_cond1}
	\end{equation}
	\item There is no step in the iteration beyond which one of the components of $\bld{x}$ stops getting updated, which we state precisely as follows.
	\begin{equation}
	    \text{Given} \,\,
	    i \in \{1,2,...,n\}\, \text{and}\, j \in \mathbb{N},\,\quad \exists\,\, l > j \,\, \text{s.t.} \,\, u(l) = i.
	    \label{eq:chaotic_cond2}
	\end{equation}
	 
\end{itemize}
In such a case the chaotic scheme is identified by the tuple $(\bld{B},\bld{C},\mathcal{S})$ where $\mathcal{S} := \{s_1, s_2, ...,s_n, u\}$.

The main result proved by Chazan and Miranker \cite{async:chazan_1969} is the following very general theorem. We state it below for completeness. Note that $\bld{v} > \bld{w}$ for two vectors $\bld{v}$ and $\bld{w}$ means that each component of $\bld{v}$ is greater than the corresponding component of $\bld{w}$, $\vert\bld{M}\vert$ is the matrix having as its entries the absolute values of the corresponding entries of $\bld{M}$, and $\rho(\bld{M})$ denotes the spectral radius of matrix $\bld{M}$.

\nomenclature{$\rho(\bld{M})$}{The spectral radius of matrix $\bld{M}$} %
\nomenclature{$\vert\bld{M}\vert$}{The matrix of the same size as $\bld{M}$ and having as its entries the absolute values of the corresponding entries of $\bld{M}$} %

\begin{theorem}
	(a) The scheme $(\bld{B},\bld{C},\mathcal{S})$ converges if $\exists \, \bld{v}\in\mathbb{R}^n$ and $\alpha < 1$ such that $\bld{v} > \bld{0}$ and $|\bld{B}|\bld{v} \leq \alpha\bld{v}$. (b) This happens if $\rho(\bld{|B|}) < 1$. (c) If no such $\bld{v}$ exists, there exists a sequence $\mathcal{S}_0$ depending on $\bld{B}$ such that the scheme $(\bld{B},\bld{C},\mathcal{S}_0)$ does not converge.
	\label{thm:chazan-miranker}
\end{theorem}
These ideas have been extended to nonlinear iterations by Baudet \cite{async:baudet_1978} and Frommer and Szyld \cite{async:frommer_2000}. These two works also present a more generalized iteration in which multiple entries are updated at each step.

We generally do not fully prescribe the shift and update functions $\mathcal{S}$ in practice, but rather let them be governed by the hardware and the non-zero pattern of the matrix $\bld{A}$. The shift and update functions could vary depending on the number of processing elements, memory access latency from different processing elements, scheduling etc. as well as the load imbalance among the updates brought about by the sparsity pattern, resulting in chaotic behaviour. The nature of the iteration could range between Jacobi and Gauss-Seidel. If shift functions are such that older and older components are used for updates until all components are updated once, the iteration becomes Jacobi, while if the shift functions are always zero, the iteration becomes Gauss-Seidel. In practice, it would be somewhere in between. This framework allows a situation in which, for the update of any component, the latest available values of all other components are used. The other components may or may not have been updated yet in the `current' solver iteration. This is useful when multiple processors are available. 

Even though we do not fully specify $\mathcal{S}$, we can influence it. One of the ways this is done is by ordering the parallel loop in a specific manner. If all the loop iterations (or work items) were actually executed in parallel, the ordering would not matter. But in general, the number of work items is much greater than the number of available processing elements and the (partial) order of execution of the work items can be influenced by the ordering of the loop.

\subsection{Asynchronous triangular solves}
\label{sec:async_triangular}

To solve a lower triangular system, the matrix $\bld{L}$ can be split as $\bld{D + \tilde{L}}$, with $\bld{D}$ diagonal non-singular and $\bld{\tilde{L}}$ strictly lower triangular. We can solve this by chaotic relaxation. A discussion of asynchronous iterations applied to triangular solves can be found in \cite{async:anzt_triangular}, in which Anzt et al. proved that an asynchronous iteration for a triangular system always converges \cite[section~2.1]{async:anzt_triangular}.

We show below the algorithm for an asynchronous lower triangular solve (algorithm \ref{alg:async-ftri}).
It may be compared with the forward substitution algorithm \ref{alg:ftrisubs}.
\begin{algorithm}[h]
	\caption{Asynchronous forward triangular solve $\bld{L}\bld{x} = \bld{b}$} 
	\label{alg:async-ftri}
	\begin{algorithmic}[1]
		\Require Lower triangular $\bld{L}$, initial solution $\bld{x}$, right-hand side vector $\bld{b}$, number of async. sweeps $n_{swp}$
		\Function{Async\_forward\_triangular\_solve}{$n_{swp} \in \mathbb{N}$, $\bld{L} \in \mathbb{R}^{n\times n}$, $\bld{b} \in \mathbb{R}^n$, $\bld{x} \in \mathbb{R}^n$}
		\State Begin parallel region and launch threads
		\For {integer $i_{swp}$ in 1..$n_{swp}$}
		\For {integer $i$ in 1..$n$} in parallel dynamically: \label{line:1iloop}
		\State $ x_i \gets (b_i - \sum_{j=1}^{i-1} L_{ij}x_j)/L_{ii} $
		\EndFor \hspace{0.2cm} (no synchronization)
		\EndFor \hspace{0.2cm} (no synchronization)
		\State End parallel region
		\State \Return $\bld{x}$
		\EndFunction
	\end{algorithmic}
\end{algorithm}
The backward triangular solve is similar, except that the loop starting at line 4 is ordered backwards.

Several `sweeps' (or `global iterations') over all unknowns are carried out, where every sweep updates each unknown once. The `dynamically parallel' loop in these algorithms implies that the work items are not divided among processing elements (cores) \emph{a priori}; rather, new work items are assigned to processing elements as and when the latter become free. Note that the loop over the sweeps is started inside the parallel region, which means each thread keeps count of its own sweeps. There is no synchronization at the end of the loop - even if some threads are computing entries for the first sweep, other threads may start executing work-items for the next sweep. We observe that for such a method to be useful, it must converge sufficiently in a small number of sweeps independent of the number of parallel processing elements.

\subsection{Asynchronous ILU factorization}
\label{sec:AILU}
Chow and Patel proposed \cite{ilu:chowpatel} a highly parallel ILU factorization algorithm based on comparing the left and right sides of the equation $[\bld{LU}]_{ij} = A_{ij}$.
Suppose we restrict the sparsity pattern of the computed incomplete LU factorization to an index set $S$, which necessarily contains the diagonal positions $(j,j) \, \forall \, j \in \{1,2,...,n\}$.  Let $m := |S|$, the number of non-zeros in the factorization. Then, the above component-wise equality leads to
\begin{equation}
\begin{split}
L_{ij} &= \left(A_{ij} - \sum_{k=1}^{j-1}L_{ik}U_{kj}\right) / U_{jj}, \qquad \text{ if } (i,j) \in S,\, i>j \\
U_{ij} &= A_{ij} - \sum_{k=1}^{i-1}L_{ik}U_{kj} \qquad \qquad \qquad \text{ if } (i,j) \in S,\, i \leq j.
\end{split}
\label{eq:scalarilu_fp_detail}
\end{equation}
This can be written as $\bld{x} = \bld{g}(\bld{x})$,
where $\bld{x} \in \mathbb{R}^m$ is a vector containing all the unknowns $L_{ij}$ and $U_{ij}$ in some order. These lead to a fixed point iteration of the form 
\begin{equation}
\bld{x}^{n+1} = \bld{g}(\bld{x}^n).
\label{eq:scalarilu_fp}
\end{equation}
An asynchronous form of the above nonlinear fixed-point iteration, as given by Chow and Patel, is shown in algorithm \ref{alg:asyncILU_chowpatel}.
\begin{algorithm}
	\caption{Asynchronous ILU factorization} \label{alg:asyncILU_chowpatel}
	\begin{algorithmic}[1]
		\Require Assign initial values to $L_{ij}$ and $U_{ij}$. 
		Let $S$ be the desired nonzero index-set.
		\For{ $i_{swp}$ in 1..$n_{swp}$}
		\For { $(i,j) \in S$} in parallel:
		\If  {$i > j$} 
		\State $L_{ij} \gets \left(A_{ij} - \sum_{k=1}^{j-1}L_{ik}U_{kj}\right) / U_{jj} $
		\Else 
		\State $ U_{ij} \gets A_{ij} - \sum_{k=1}^{i-1}L_{ik}U_{kj} $
		\EndIf
		\EndFor
		\EndFor
	\end{algorithmic}
\end{algorithm}
Note that the ordering of the inner parallel loop still matters when the number of processing elements is less than the number of non-zeros in $S$. This ordering can be chosen for more effective preconditioning rather than to alleviate data dependence. Depending on the local imbalance and number of processors, the resulting iteration could range from nonlinear Jacobi to nonlinear Gauss-Seidel. If one processor is used (resulting in nonlinear Gauss-Seidel) and the loop is ordered in `Gaussian elimination' ordering (eg. row-wise or column-wise ordering), traditional ILU is recovered. As an initial guess, we use the entries of $\bld{A}$. In a CFD simulation, the approximate LU factorization from the previous time step can also be used as the initial guess.

We briefly describe the framework used by Chow and Patel. Let an ordering of the unknowns $\bld{x}$ (the entries of the triangular factors) be given by the bijective map $\alpha:S\rightarrow \{1,2,3,...,m\}$ (where $m=|S|$). $\bld{x}$ can now be expressed as
\begin{equation}
x_{\alpha(i,j)} = 
\begin{cases}
L_{ij} \quad &\text{ if } i > j \\
U_{ij} \quad &\text{ if } i \leq j
\end{cases}.
\label{eq:scalarordering}
\end{equation}

\nomenclature{$S$}{Sparsity pattern made up of indices $(i,j)$ of non-zero entries}

Further, the mapping  $\bld{g}:D\rightarrow\mathbb{R}^m$ can be expressed, for $(i,j) \in S$, as
\begin{equation}
g_{\alpha(i,j)}(\bld{x}) = \begin{cases}
(A_{ij} - \sum_{k = 1}^{j-1} x_{\alpha(i,k)}x_{\alpha(k,j)})/x_{\alpha(j,j)} \quad &\text{if } i > j \\
A_{ij} -  \sum_{k = 1}^{i-1} x_{\alpha(i,k)}x_{\alpha(k,j)}  &\text{if } i \leq j
\end{cases}.
\end{equation}
The domain of definition of $\bld{g}$ is 
$D := \{\bld{x}\in\mathbb{R}^m \,|\, x_{\alpha(j,j)} \neq 0 \,\forall \,j \in \{1,2,...,m\} \}$.

We now select any one of the `Gaussian elimination orderings' for $\alpha$, such as
\begin{equation}
(1,1) \prec (1,2) \prec ... \prec (1,n) \prec (2,1) \prec (2,2) \prec ... \prec (n,n-1) \prec (n,n).
\label{eq:geordering}
\end{equation}
This is a row-major ordering; other possible Gaussian elimination orderings are column-major ordering and the partial ordering chosen by Chow and Patel \cite{ilu:chowpatel}.

Chow and Patel proved \cite[theorem~3.5]{ilu:chowpatel} local convergence of the asynchronous ILU iteration; that is, a fixed point of $\bld{g}$ is a point of attraction of the iteration.
In the event of asynchronous updates, some diagonal entries may sometimes be set to zero. To take this into account, Chow and Patel define a `modified Jacobi-type iteration' corresponding to the iteration in equation \eqref{eq:scalarilu_fp}, in which whenever a zero diagonal entry is encountered, it is replaced by an arbitrary non-zero value.
They showed \cite[theorem~3.7]{ilu:chowpatel} that a synchronized modified Jacobi-type iteration corresponding to equation \eqref{eq:scalarilu_fp} converges in at most $m$ iterations irrespective of the initial guess.

\begin{algorithm}
	\caption{Asynchronous ILU factorization} \label{alg:asyncILU}
	\begin{algorithmic}[1]
		\Require Assign initial values to $L_{ij}$ and $U_{ij}$. 
		Let $S$ be a set of non-zero $i,j$ indices.
		\State Begin parallel region
		\For{ $i_{swp}$ in 1..$n_{swp}$}
		\For { $i$ in 1..$n$} in parallel dynamically:
		\For {$j$ in 1..$n$, $(i,j) \in S$}
		\If  {$i > j$} 
		\State $L_{ij} \gets \left(A_{ij} - \sum_{k=1}^{j-1}L_{ik}U_{kj}\right) / U_{jj} $ \label{line:L}
		\Else 
		\State $ U_{ij} \gets A_{ij} - \sum_{k=1}^{i-1}L_{ik}U_{kj} $ \label{line:U}
		\EndIf
		\EndFor
		\EndFor \hspace{0.2cm} (no synchronization)
		\EndFor \hspace{0.2cm} (no synchronization)
		\State End parallel region
	\end{algorithmic}
\end{algorithm}

In our implementation (algorithm \ref{alg:asyncILU}), the loop over entries of the matrix is always ordered in the row-major ordering. Each work-item consists of computing all the entries in one row of the matrix. 
In OpenMP, a number of consecutively ordered work-items (determined by the loop ordering) are placed in one `chunk'. In case of dynamic scheduling of work, when a free thread requests work, it is assigned the next chunk of waiting work-items \cite[section~2.7.1]{openmp45}. For our problem, each chunk consists of a number of consecutive rows of the matrix.
Thus, the row-major ordering of the loop means that all entries in the rows in one chunk are computed in a Gaussian elimination ordering (the row-major ordering) with respect to entries in that chunk.

We point out that no atomic operations are used in our implementation. The aim was to write code as close as possible to the mathematical notion of asynchronous iterations. It is not clear whether this requires atomic operations. Atomic writes could be used for lines \ref{line:L} and \ref{line:U} in algorithm \ref{alg:asyncILU}, and this would ensure ``valid" numbers (avoiding corruption) after each step. This amounts to imposing a local synchronization between two threads, one of which is attempting to write to a location that the other is attempting to read. This has not been investigated in this work. Using such atomic operations may affect the results presented in section \ref{sec:results}.
A related issue, which is not evident in algorithm \ref{alg:asyncILU}, is that of accumulating the sums in lines \ref{line:L} and \ref{line:U}. We found that accumulating the sum directly in the memory location of $L_{ij}$ caused much slower convergence of the asynchronous iteration. However, if a local variable accumulated the sum first, and the final value of this variable was assigned to $L_{ij}$, convergence is faster and more reliable. This may be because in the former approach, there are times when $L_{ij}$ contains a partial sum, and is read by another thread for the update of another component. A partial sum can be very different from a stale (older) value of $L_{ij}$, and does not conform to the definition \ref{eq:chaotic_relaxation} of asynchronous iteration. Thus, two pieces of code that represent completely equivalent sequential algorithms may yield different asynchronous iterations.

To apply the ILU preconditioner (ie., to solve $\bld{LUx}=\bld{b}$ given $\bld{L}$ and $\bld{U}$) in parallel, Chow and Patel used synchronous Jacobi iterations to solve the triangular systems $\bld{L}\bld{y} = \bld{b}$ and $\bld{U}\bld{x} = \bld{y}$ \cite{ilu:chowpatel}.
\section{Asynchronous block preconditioners}
\label{sec:async_block}
As we discussed in section \ref{sec:solver}, our Jacobian matrix has a natural block structure with small dense blocks. We would like to take advantage of this fact in asynchronous preconditioning. As we show later, asynchronous block iterations are generally more robust and converge in fewer iterations than their scalar counterparts. Further, blocking may improve cache utilization and vectorization.

\subsection{Asynchronous block triangular solves}
We can extend the asynchronous triangular solves to block triangular solves in a straightforward manner. We present the forward block triangular solve algorithm as an example. Assume that the blocks are square $b\times b$ and that the matrix size $n$ is divisible by the block size $b$. $\bld{x}_{a:b}$ denotes the subvector $[x_a, x_{a+1}, ..., x_{b-1}]^T$ of $\bld{x}$ for $a<b$, so that $\bld{x}_{(i-1)b+1:ib+1}$ denotes the $i$-th subvector of length $b$. Similarly, $\bld{L}_{(i-1)b+1:ib+1,\,(j-1)b+1:jb+1}$ denotes the $(i,j)$th $b\times b$ sub-block of the matrix $\bld{L}$. Let us simplify this cumbersome notation by defining $\bld{x}_i := \bld{x}_{(i-1)b+1:ib+1}$ and $\bld{L}_{ij} := \bld{L}_{(i-1)b+1:ib+1,\,\,(j-1)b+1:jb+1}$. 
\begin{algorithm}[H]
	\caption{Asynchronous forward block-triangular solve $\bld{L}\bld{x} = \bld{b}$} 
	\label{alg:b-async-ftri}
	\begin{algorithmic}[1]
		\Require $\bld{L}$ is lower block triangular with nonsingular diagonal blocks, $\bld{x}$ is an initial guess for the solution, $\bld{b}$ is the right-hand side vector, $b$ is the size of square blocks in $\bld{L}$
		\Function{Async\_forward\_block\_triangular\_solve}{$n_{swp} \in \mathbb{N}$, $\bld{L} \in \mathbb{R}^{n\times n}$, $\bld{b} \in \mathbb{R}^n$, $\bld{x} \in \mathbb{R}^n$}
		\State Begin parallel region and launch threads
		\For {integer $i_{swp}$ in 1..$n_{swp}$}
		\For {integer $i$ in $1$..$n/b$} in parallel dynamically:
		\State $\bld{x}_i \gets \bld{L}_{ii}^{-1}(\bld{b}_i -\sum_{j=1}^{i-1} \bld{L}_{ij}\bld{x}_j)$
		\EndFor \hspace{0.25cm} no synchronization
		\EndFor
		\State End parallel region
		\State \Return $\bld{x}$
		\EndFunction
	\end{algorithmic}
\end{algorithm}

Following the discussion on (scalar) triangular solves in \cite[section~2.1]{async:anzt_triangular}, we can easily extend the proof of convergence to this asynchronous block triangular iteration.
%

\begin{theorem}
The chaotic relaxation for the system $\bld{Lx}=\bld{b}$ defined by the splitting $\bld{D + \tilde{L}}$, with $\bld{D}$ nonsingular block diagonal and $\bld{\tilde{L}}$ strictly lower block triangular, converges.
\label{thm:block_async_triangular_local}
\end{theorem}
\begin{proof}
The Jacobi iteration for the splitting can be written as
\begin{equation}
\bld{x}^{j+1} = \bld{D}^{-1}\bld{b} - \bld{D}^{-1}\bld{\tilde{L}}\bld{x}^j.
\label{eq:blocktri}
\end{equation}
Since $\bld{\tilde{L}}$ is strictly lower block triangular and $\bld{D}^{-1}$ is block diagonal, $\bld{D}^{-1}\tilde{\bld{L}}$ is also strictly lower block triangular.
Thus, $\rho(|\bld{D}^{-1}\tilde{\bld{L}}|) = \rho(\bld{D}^{-1}\tilde{\bld{L}}) = 0$. By theorem \ref{thm:chazan-miranker}, the chaotic relaxation \eqref{eq:chaotic_relaxation} for this splitting converges.
\end{proof}

\begin{theorem}
The Jacobi iteration \eqref{eq:blocktri} for the lower block triangular linear system $\bld{Lx}=\bld{b}$ converges to the solution $\bld{L}^{-1}\bld{b}$ in at most $n/b$ iterations irrespective of the initial guess and order of updates, where $n$ is the number of unknowns in the system and $b$ is the block size.
\end{theorem}
\begin{proof}
We proceed to prove this theorem by induction. Since the iteration matrix $\bld{D}^{-1}\tilde{\bld{L}}$ is strictly block lower triangular, the $b$ unknowns in the first block $\bld{x}_1$ do not depend on any other entries of $\bld{x}$.
Thus, irrespective of the order of updates and the initial guess, after the first iteration of equation \eqref{eq:blocktri}, the unknowns in the first block attain their exact values $\bld{x}_1^1 = (\bld{D}^{-1}\bld{b})_1$.

Let us assume that the first $j-1$ blocks' unknowns have attained their exact values by the end of the $(j-1)$th iteration. The unknowns of the $j$th block $\bld{x}_j$ depend only on those of blocks 1,2,...,$j-1$ due to the strictly block lower triangular nature of the iteration matrix. 
Therefore, the $j$th block attains its exact values by the end of the $j$th iteration.

Thus, by induction, all $n/b$ blocks' unknowns attain their exact values by the end of the $(n/b)$th iteration.
\end{proof}

We now define a `block-asynchronous' iteration. Let the vector function $\bld{\psi}:D\subset\mathbb{R}^n \rightarrow\mathbb{R}^n$ be any iteration. For the purpose of forward triangular solves, it is defined by $\bld{\psi}(\bld{x}) := \bld{D}^{-1}\bld{b} - \bld{D}^{-1}\tilde{\bld{L}}\bld{x}$. The definition is inspired by Frommer and Szyld's definition of asynchronous iteration \cite[definition~2.2]{async:frommer_2000}.
Again, $\psi_i$ denotes the $i$th function while $\bld{\psi}_i$ (in bold face) denotes the sub-vector of $b$ functions in the $i$th block.
\begin{definition} A block-asynchronous iteration is of the form
\begin{equation}
\bld{x}_i^{j+1} = 
\begin{cases} 
\bld{x}_i^j & \text{if } i \neq u(j) \\
\bld{\psi}_i(x_1^{j-s_1(j)}, x_2^{j-s_2(j)},...,x_n^{j-s_n(j)}) & \text{if } i = u(j)
\end{cases}
\label{eq:chaotic_blocklowertri}
\end{equation}
where $s_\alpha:\mathbb{N}\rightarrow\mathbb{N},\, \alpha \in \{1,2,...,n\}$ are the shift functions and $u:\mathbb{N}\rightarrow\{1,2,...,n/b\}$ is the update function. ($n/b$ is an integer equal to the number of cells.)
\end{definition}
As before, the shift functions are assumed to satisfy criterion \eqref{eq:chaotic_cond1}, while the update function satisfies a block version of \eqref{eq:chaotic_cond2}:
\begin{equation}
	    \text{Given} \,\,
	    i \in \{1,2,...,n/b\}\, \text{and}\, j \in \mathbb{N},\,\quad \exists\,\, l > j \,\, \text{s.t.} \,\, u(l) = i.
	    \label{eq:blockchaotic_cond2}
\end{equation}
Recall that a step is defined by a single read of (all the required entries from) the unknown vector $\bld{x}$ by one thread. Here, this is followed by the update of a point-block of entries by the thread. Hence, the range of the update function consists of integers between 1 and the number of blocks. This allows exact treatment of dependencies of variables within the point-block.
However, the individual entries used for computing the iterate can come from different previous updates. Therefore there is a shift function for every entry, instead of every block, of the unknown vector.
This block-asynchronous iteration is equivalent to Frommer and Szyld's asynchronous iteration \cite[definition~2.2]{async:frommer_2000} in the case when the set of indices updated in each step exactly corresponds to one point-block. As such, the established asynchronous convergence theory given by Frommer and Szyld \cite{async:frommer_2000} still holds.

\begin{theorem} The block-asynchronous linear iteration \eqref{eq:chaotic_blocklowertri} for solving the lower block triangular system converges in a finite number of steps to the solution $\bld{L}^{-1}\bld{b}$.
\label{thm:blockasync_triangular_global}
\end{theorem}
\begin{proof}
Because of the strictly lower block-triangular nature of the iteration matrix, the first cell's unknowns (the $b$-block $\bld{x}_1$) do not depend on any other unknowns. Whenever this block is updated for the first time, it attains its exact values.

Next, let us assume that at the end of the $k\textsuperscript{th}$ step, all cells up to the $j\textsuperscript{th}$ cell have attained their exact values. Since $\hat{s}$ is the upper bound on the delays (by condition \eqref{eq:chaotic_cond1}), starting latest at the $k+\hat{s}+1$ step, the exact values for all cells up to the $j\textsuperscript{th}$ cell will be used for any update that requires them.
By condition \eqref{eq:chaotic_cond2}, we know that there exists some step $l \geq k+\hat{s}+1$ for which $u(l) = j+1$. That is, there exists some step $l$ at which the unknowns of the $j+1\textsuperscript{th}$ cell's block is updated using the exact values for all the cells that it depends on. Thus the $j+1\textsuperscript{th}$ cell attains its exact values after a finite number of steps.

Hence, by induction, all the unknowns equal their fixed-point values after some finite number of steps.
\end{proof}

Even though asynchronous iteration or Chazan-Miranker's chaotic relaxation does not have global iterations, our implementation is carried out using several `sweeps' of a parallel OpenMP loop over all the unknowns. If we use $n_s$ sweeps, every unknown is updated exactly $n_s$ times.

Since $\bld{Ly}=\bld{b}$ and $\bld{Ux}=\bld{y}$ need to be solved every time the preconditioner is applied in a solver, the above bounds on the number of fixed-point iterations or steps are only of academic interest. We hope to be able to use only a few sweeps of asynchronous updates to approximately solve the triangular systems.

It must be mentioned that Anzt et al. \cite{async:anzt_blockasync} used the term `block-asynchronous iteration' to describe a different kind of iteration that they proposed. In their approach to linear asynchronous iterations, the block aspect accounts for features of the hardware or the programming model (thread blocks in CUDA, in their case). They perform several Jacobi iterations within a block to invert them approximately, while the coupling between the blocks has asynchronous character. In our case, a diagonal block is inverted exactly and sequentially by one CPU thread, because our blocks are smaller and dense.
The above block-asynchronous iteration is also applicable to non-linear iterations, as required by the iterative block-ILU method. The similarity between our framework and the iteration of Anzt et al. is that the coupling between blocks is asynchronous.

\subsection{Asynchronous block ILU preconditioner}

Using the equation $(\bld{LU})_{ij} = \bld{A}_{ij}$, we can define a block LU factorization. Suppose $S_B$ is a set of block indices which includes all diagonal blocks; it will be the block sparsity pattern imposed on the computed block LU factors. The block ILU factorization is computed as
\begin{equation}
\begin{split}
\bld{L}_{ij} &= \left(\bld{A}_{ij} - \sum_{k=1}^{j-1}\bld{L}_{ik}\bld{U}_{kj}\right) \bld{U}_{jj}^{-1}, \qquad \text{ if } (i,j) \in S_B,\, i>j \\
\bld{U}_{ij} &= \bld{A}_{ij} - \sum_{k=1}^{i-1}\bld{L}_{ik}\bld{U}_{kj} \qquad \quad \qquad \quad \text{ if } (i,j) \in S_B,\, i \leq j
\end{split}
\label{eq:bilu}
\end{equation}
where subscripts denote block-indices. 

\nomenclature{$S_B$}{Block sparsity pattern made up of indices $(i,j)$ of non-zero blocks}



The unknowns are
$[\bld{L}_{ij}]_{kl},\, (i,j) \in S_B\, \text{s.t.}\, i>j$ and $[\bld{U}_{ij}]_{kl},\, (i,j) \in S_B\, \text{s.t.}\, i \leq j$, for $1 \leq k,l \leq b$.
Here, $[\bld{L}_{ij}]_{kl}$ denotes the $(k,l)$ entry of the sub-matrix $\bld{L}_{ij}$.
For analyzing block preconditioners, we introduce an ordering of the unknowns by the bijective map
\begin{equation}
\beta:S_B \times \{1,2,...,b\} \times \{1,2,...,b\} \rightarrow \{1,2,3,...,m\},
\end{equation}
where $m=|S_B|b^2$, the total number of non-zeros in the block-ILU factorization. We can then write a vector of all the unknowns and call it $\bld{x}$, which is ordered as (compare with $\alpha$ defined in \eqref{eq:scalarordering})
\begin{equation}
x_{\beta(i,j,k,l)} = 
\begin{cases}
[\bld{L}_{ij}]_{kl} \quad & \text{ if } i > j \\
[\bld{U}_{ij}]_{kl} \quad & \text{ if } i \leq j
\end{cases}.
\end{equation}
We can also define a block-ordering $\beta_B:S_B \rightarrow \{1,2,...,m/b^2\}$,
which is bijective and related to the ordering $\beta$ by $\beta_B(i,j) = \beta(i,j,b,b)/b^2 $.

Let us define $\bld{h}:D_B\rightarrow\mathbb{R}^m$ (where $D_B \subset \mathbb{R}^m$) to be the function that represents the right-hand-side of equation \eqref{eq:bilu}.
Denote the matrix $\bld{x}_{\beta(i,j,:,:)}$ by $\bld{X}_{ij} \in \mathbb{R}^{b\times b}$ for $(i,j) \in S_B$. Note that `:' at an indexing position denotes the entire range of indices possible for that position; here it denotes all integers from 1 to $b$. Similarly, we denote $\bld{h}_{\beta(i,j,:,:)}:D_B \rightarrow \mathbb{R}^{b\times b}$ by $\bld{H}_{ij}$ for $(i,j) \in S_B$.
Then we can express the domain of definition of $\bld{h}$ as
\begin{equation}
D_B := \{\bld{x}\in\mathbb{R}^m \,|\, \bld{X}_{jj} \text{ is nonsingular } \,\forall \,j \in \{1,2,...,m/b^2\} \}.
\end{equation}
With this, the mapping  $\bld{h}:D_B\rightarrow\mathbb{R}^m$ can be expressed as
\begin{equation}
\bld{H}_{ij}(\bld{x}) := \bld{h}_{\beta(i,j,:,:)}(\bld{x}) = \begin{cases}
(\bld{A}_{ij} - \sum_{k = 1}^{j-1} \bld{X}_{ik}\bld{X}_{kj})\bld{X}_{jj}^{-1} \quad \text{if } i > j \\
\bld{A}_{ij} -  \sum_{k = 1}^{i-1} \bld{X}_{ik}\bld{X}_{kj}  \qquad \text{if } i \leq j
\end{cases}.
\label{eq:blockH}
\end{equation}
Now the fixed-point of the ILU iteration can be written as $\bld{x} = \bld{h}(\bld{x})$,
the synchronized Jacobi-type fixed-point iteration can be expressed as
\begin{equation}
\bld{x}^{n+1} = \bld{h}(\bld{x}^n),
\label{eq:seq_blockilu_fp}
\end{equation}
while the block-asynchronous iteration is expressed as
\begin{equation}
\bld{X}_{ij}^{k+1} =
    \begin{cases}
        \bld{H}_{ij}(x_1^{k-s_1(k)}, x_2^{k-s_2(k)}, ..., x_n^{k-s_m(k)}), \quad & \beta_B(i,j) = u(k) \\
        \bld{X}_{ij}^k \quad \quad & \beta_B(i,j) \neq u(k)
    \end{cases},
\label{eq:async_blockilu_fp}
\end{equation}
where 
$s_j:\mathbb{N}\rightarrow\mathbb{N}$ (for $1\leq j \leq m$) are the shift functions and $u:\mathbb{N}\rightarrow\{1,2,...,m/b^2\}$ is the update function. The shift and update functions are assumed to have the following properties, corresponding to the usual properties of asynchronous iterations.
 \begin{equation}
	    \exists\, \hat{s} \in \mathbb{N}\,\, \text{s.t.}\,\, 0 \leq s_i(k) \leq \min\{k-1,\hat{s}\} \quad \forall i \in \{1,2,...m\},\, k \in \mathbb{N}.
	    \label{eq:blockchaotic_ilu_cond1}
\end{equation}
\begin{equation}
	    \text{Given} \,\,
	    i \in \{1,2,...,m/b^2\}\, \text{and}\, k \in \mathbb{N},\,\, \exists\,\, l > k \,\, \text{s.t.} \,\, u(l) = i.
	    \label{eq:blockchaotic_ilu_cond2}
\end{equation}

\nomenclature{$\mathbb{N}$}{The set of all natural numbers 1,2,...}

Let us select a `block Gaussian elimination' ordering for $\beta_B$ as the same as the Gaussian elimination ordering in \eqref{eq:geordering}, except that the indices now correspond to block indices.
\begin{equation}
(1,1) \prec (1,2) \prec ... \prec (1,\frac nb) \prec (2,1) \prec (2,2) \prec ... \prec (\frac nb,\frac nb-1) \prec (\frac nb, \frac nb), \quad (i,j) \in S_B
\label{eq:bgeordering}
\end{equation}
($n$ is the dimension of the matrix $\bld{A}$). This defines a partial ordering for $\beta$, where the ordering within each block is unspecified. We can now extend the convergence analysis given by Chow and Patel for the asynchronous ILU process to the asynchronous block-ILU algorithm.

\begin{lemma}
$\bld{h}(\bld{x})$ is differentiable for all $\bld{x}\in D_B$.
\label{lma:differentiable}
\end{lemma}
\begin{proof}
From the definition \eqref{eq:blockH} we see that each block $\bld{H}_{ij}$ is a rational matrix function of the $\bld{X}_{kl}$. This implies that each $h_i$ is a rational function of the $x_j$.
Thus, $\bld{h} \in [C^1(D_B)]^m$, that is, $\bld{h}$ is continuously differentiable in its domain of definition.
\end{proof}
\begin{lemma}
When the unknowns $\bld{x}$ and the mappings $\bld{h}$ are ordered in the block Gaussian elimination ordering \eqref{eq:bgeordering}, $\bld{h}$ has strictly lower-triangular structure, ie., $h_k(\bld{x})$ depends only on $\{x_1,x_2,...,x_{k-1}\}$.
Thus the Jacobian $\bld{h}'(\bld{x})$ is strictly lower triangular $\forall \,\bld{x}\in D_B$.
\label{lma:blowertriangular}
\end{lemma}
\begin{proof}
From equation \eqref{eq:blockH}, we see that $\bld{H}_{ij}$ depends on $\bld{X}_{pq}$ only if 
\begin{equation}
(p,q) \in \{(\gamma,\delta) \in S_B |\, \gamma < i,\, \delta = j\} \cup \{(\gamma,\delta) \in S_B | \, \gamma = i,\, \delta < j\}.
\end{equation}
This implies that $\bld{H}_{ij}$ depends on a subset of the blocks $\{\bld{X}_{pq} \,|\, p \leq i, q \leq j,\, (p,q) \neq (i,j)\}$ preceding it in the ordering \eqref{eq:bgeordering}.
This means that for indices in the $(i,j)$th block, ie., for indices $k$ such that $\beta(i,j,1,1) \leq k \leq \beta(i,j,b,b)$, $h_k(\bld{x})$ only depends on a subset of $\{x_1,x_2,...x_{b^2 [(k-1)/b^2]}\} \subseteq \{x_1,x_2,...x_{k-1}\}$. Here, $[.]$ is the greatest integer function.
Thus $\bld{h}$ has not only a strictly lower triangular structure but a strictly lower block triangular structure.

Therefore, each $h_k$ is a function of $x_j$ only for $j<k$ :
\begin{equation}
h_k(\bld{x}) = h_k(x_1, x_2,...,x_{k-1}) \implies \frac{\partial h_k}{\partial x_l} = 0 \quad \forall\, l \geq k,
\end{equation}
which means the Jacobian $\bld{h}'(\bld{x})$ is strictly lower triangular $\forall \, \bld{x} \in D_B$.
\end{proof}

\begin{theorem}
The synchronous nonlinear fixed-point iteration $\bld{x}^{p+1} = \bld{h}(\bld{x}^p)$ is locally convergent, that is, if $\bld{x}_*$ is a fixed point of $\bld{h}$, it is a point of attraction of the synchronous iteration.
\end{theorem}
\begin{proof}
By lemmas \ref{lma:differentiable} and \ref{lma:blowertriangular}, $\bld{h}$ is differentiable and all eigenvalues of $\bld{h}'$ are zero. Since eigenvalues are unaffected by symmetric permutations of the matrix, this holds for any reordering of the equations and unknowns.

Thus, the spectral radius of $\bld{h}'$ is zero. Now from the Ostrowski theorem (theorem 10.1.3 of Ortega and Rheinboldt \cite{ortega_rheinboldt}),if $\bld{x}_*$ is a fixed point of $\bld{h}$, it is a point of attraction.
\end{proof}

\begin{theorem}
The asynchronous iteration corresponding to \eqref{eq:seq_blockilu_fp} is locally convergent, ie., if $\bld{x}_*$ is a fixed point of the asynchronous iteration, it is a point of attraction of the iteration.
\label{thm:async_blockilu_localconv}
\end{theorem}
\begin{proof}
In lemmas \ref{lma:differentiable} and \ref{lma:blowertriangular}, we showed that $\bld{h}$ is differentiable and all eigenvalues of $\bld{h}'$ are zero. Therefore, $\rho(|\bld{h}'(\bld{x}_*)|) = 0$.
Note that for a matrix $\bld{A}$, $|\bld{A}|$ denotes the matrix of absolute values of the corresponding entries.

Thus, according to the theorem \cite[theorem 4.4]{async:frommer_2000} by Frommer and Szyld, $\bld{x}_*$ is a point of attraction of the asynchronous iteration corresponding to \eqref{eq:seq_blockilu_fp}.
\end{proof}

Similar to Chow and Patels's work \cite{ilu:chowpatel}, global convergence of the synchronized block ILU \eqref{eq:seq_blockilu_fp} can be proved.

\begin{theorem}
If a fixed point of the function $\bld{h}$ exists, it is unique.
\end{theorem}
\begin{proof}
We assume the Gaussian elimination ordering \eqref{eq:bgeordering}. This leads to no loss of generality because $\bld{x}_*$ is a fixed point of $\bld{h}$ in one ordering if and only if it is a fixed point in another ordering.

In this ordering, the solution can be found in one iteration through forward substitution due to the strictly lower triangular nature of $\bld{h}$. The forward substitution completes if no $\bld{X}_{jj} = \bld{U}_{jj}$ is set to a singular matrix, in which case the solution is the unique fixed point because forward substitution gives a unique solution. If any of these diagonal blocks is set to a singular matrix, no fixed point exists.
\end{proof}

In analogy with the scalar ILU theory, a modified iteration corresponding to the sequential or asynchronous iteration can be defined as a similar iteration except that when a diagonal block $\bld{X}_{jj}$ becomes singular, it is replaced by an arbitrary non-singular matrix.

\begin{theorem}
If $\bld{h}$ has a fixed point, the modified Jacobi-type iteration corresponding to \eqref{eq:seq_blockilu_fp} converges in at most $m/b^2$ iterations from any initial guess $\bld{x}^0$. ($m$ is the number of nonzero entries in all blocks in the sparsity pattern $S_B$ and $b$ is the block size.)
\label{thm:blockilu_globalconv}
\end{theorem}
\begin{proof}
It can be observed from \eqref{eq:blockH} that the first block $\bld{X}_{1,1}^1$ or $\bld{U}_{1,1}^1 \,:= \bld{x}_{\beta(1,1,:,:)}^1 $ of the first iterate does not depend on $\bld{x}$ at all. Thus the $b^2$ entries in $\bld{X}_{1,1}$ attain their exact fixed-point values after the first iteration.
Other blocks may potentially be updated in some arbitrary manner, though the modified iteration ensures that diagonal blocks remain non-singular and the iteration does not break down.

Suppose the next $b^2$ entries ($x_{b^2}$ to $x_{2b^2}$) in the ordering \eqref{eq:bgeordering} correspond to the block $(k,l)\in S_B$. Then $\bld{X}_{k,l}$ depends only on $\bld{X}_{1,1}$ (at most), therefore its entries attain their fixed-point values after iteration 2, and retain those values after further iterations. Thus, the first $2b^2$ entries of $\bld{x}$ attain their correct fixed-point values at the completion of iteration 2.

If all blocks up to the $p\textsuperscript{th}$ non-zero block ($p\geq 1)$ in the Gaussian elimination ordering \eqref{eq:bgeordering} have attained their final values by the  $p\textsuperscript{th}$ iteration, the $p+1\textsuperscript{th}$ non-zero block in that ordering depends only on the known blocks, due to the strictly lower triangular nature of $\bld{h}$. Thus, the  $p+1\textsuperscript{th}$ non-zero block attains its exact value by the $p+1\textsuperscript{th}$ iteration.

Continuing in this manner to the last block in the ordering \eqref{eq:bgeordering}, we conclude by induction that
all $m$ entries of $\bld{x}$ attain their fixed-point values at the completion of iteration $m/b^2$.
\end{proof}
We note that the maximum number of global iterations needed for modified Jacobi-type block ILU is smaller than that for the original modified Jacobi-type scalar ILU iteration.

For the asynchronous modified iteration, we cannot make claims about the number of global iterations it takes for convergence because there are no global iterations. However, because of the property \eqref{eq:chaotic_cond2} of asynchronous iteration, for each block there are steps at which it continues to get updated as we perform more and more asynchronous steps. By an argument similar to the proof of theorem \ref{thm:blockasync_triangular_global}, we can conclude that the asynchronous modified iteration will converge in a finite number of steps.


To parallelize the application of the asynchronous ILU preconditioner, we use asynchronous triangular solves as shown in algorithm \ref{alg:async-ftri}, and for the asynchronous block ILU preconditioner we use asynchronous block-triangular solves as illustrated in algorithm \ref{alg:b-async-ftri}.

As mentioned earlier, OpenMP divides the work items into chunks. On a CPU, whenever an idle thread is assigned work, it is assigned an entire chunk which is then processed sequentially by the thread. Therefore an interesting detail to note is that we can expect a parallel (modified) fixed-point iteration to converge in at most the same number of sweeps as the number of chunks, as long as the iteration function is strictly lower triangular in the loop ordering.
The reasoning is that once the entries in all chunks before the $p\textsuperscript{th}$ chunk attain their final values, all entries in the $p+1\textsuperscript{th}$ chunk attain their final values in the next sweep.
The number of chunks is usually much less than the number of non-zero blocks, because a chunk usually contains a substantial number of work-items for performance reasons.
\section{Orderings of mesh cells}

As we will see in the results section, asynchronous (block) ILU is quite sensitive to the ordering of mesh cells. We have used the following topological orderings:
\begin{itemize}
    \item Reverse Cuthill-McKee (RCM) \cite{orderings:rcm, orderings:george_liu_book}. This is a common ordering used to solve PDEs with ILU preconditioning; it is a `level-set' ordering that aims to reduce the bandwidth of the matrix.
    \item One-way dissection (1WD) \cite{orderings:1wd_unstructured}. This algorithm aims to order the grid by recursively introducing separators and sub-dividing the grid.
\end{itemize}
The implementations available in PETSc \cite{petsc-web-page} were used to achieve these orderings. For efficiency, we reorder the grid itself in a pre-processing stage and avoid reordering of matrices during the non-linear solve.

In addition to the above-mentioned algorithms, we also introduce the following algorithms which are suited to viscous fluid dynamic simulations.

\subsection{Line ordering}
Line solvers are well established in CFD \cite{line:dplr_1998, line:incompressible, multigrid:unstructured}. Meshes for viscous flows are usually generated with high grid-stretching in the boundary layer near the body being studied (figure \ref{fig:stretched_mesh}), so as to capture the highly anisotropic flow profiles in that region efficiently.
Thus, the wall-normal direction, having a high density of points, is one of strong coupling, while the the wall-tangent direction is more loosely coupled. 
Mavriplis \cite{multigrid:unstructured} used an algorithm for finding lines of strong coupling in the grid, based on the physical locations of the grid points.
In that work, lines made up of grid vertices are found for a vertex-centred discretization. Since we use a cell-centred discretization, we apply the same algorithm to cell-centres to find lines of cells which are tightly coupled. Note that only cells that lie in regions of high anisotropy are incorporated into lines. Thus lines start at boundary cells and continue towards the interior only while a local anisotropy threshold is met. These truncated lines are sometimes called linelets in literature (eg. \cite{line:incompressible}). For simplicity, unlike Mavriplis, we only consider lines near boundaries and do not consider shear layers in the interior of the flow; our lines are limited to the boundary layer region. In this way, the grid is divided into lines of anisotropic cells and individual isotropic cells.

Once lines are found, the grid is reordered so that cells that make up a line are contiguous in the ordering. The Jacobian matrix can then be broken up into line blocks (which are block tridiagonal, representing the one-dimensional strong coupling along the line), point (cell) blocks for cells not belonging to any line, and the blocks coupling them. In a traditional line solver, a block-Jacobi iteration is applied. The individual line blocks are inverted exactly using Thomas' algorithm, and point blocks are also easily inverted exactly, but the couplings between neighbouring lines and isotropic points are not part of the preconditioner.

In our code, we reorder the grid such that cells belonging to a line are contiguous, and then apply sequential or asynchronous ILU preconditioning to this reordering. Note that a sequential ILU preconditioner would invert the line-blocks exactly. Moreover, it would more accurately relax the coupling between the lines and isotropic cells, compared to a traditional line solver.

\subsection{Hybrid line-X ordering}
\begin{figure}[h]
\begin{minipage}{0.29\linewidth}
        \includegraphics[scale=0.15]{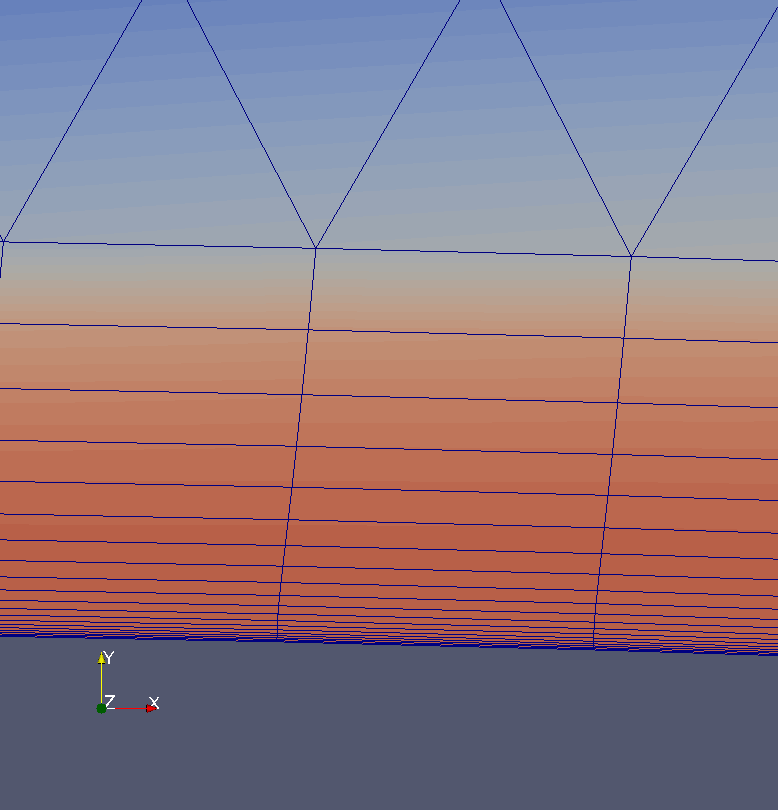}
        \captionof{figure}{Anisotropic cells \\near the boundary}
        \label{fig:stretched_mesh}
\end{minipage}
\begin{minipage}{0.71\linewidth}
        \centering
        \includegraphics[scale=0.25]{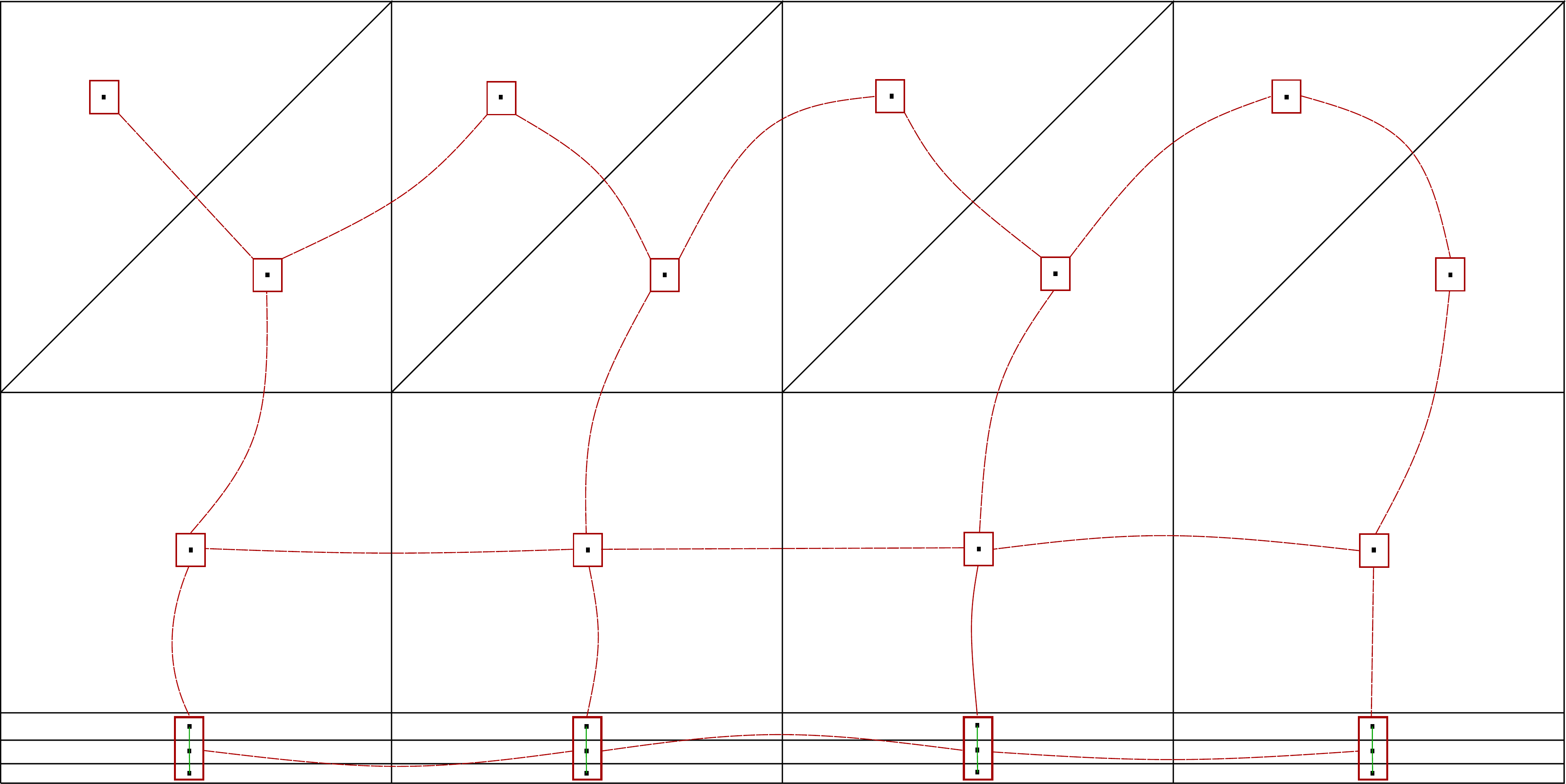}
        \captionof{figure}{The graph used to generate hybrid line-X orderings. Black lines represented mesh edges, and red boxes represent vertices of the graph and red curves represent edges of the graph. Notice that cells in the bottom three layers of anisotropic cells are connected vertically to form lines.}
        \label{fig:hybrid_graph}
\end{minipage}
\end{figure}

The line ordering described above can be combined with a topological ordering to produce potentially even better orderings. To do this, we first find lines of anisotropic cells and order the mesh according to the line ordering described above. We then define a graph whose vertices are the lines and the individual isotropic cells; that is, each graph vertex represents either a line or an isotropic cell.
If for a pair of lines, one of the lines contains a cell that neighbours at least one cell in the other line, those two lines are assumed connected directly in the graph. Similarly, isotropic cells which neighbour a cell belonging to a line are considered connected to that line in the graph. Two isotropic cells which are neighbours in the original grid are also connected in the graph (figure \ref{fig:hybrid_graph}).

This graph is then ordered by a topological ordering such as RCM or 1WD, to produce line-RCM or line-1WD ordering respectively. Note that cells in a line remain contiguous in the final ordering - thus hybrid line-X orderings are also line orderings. We will see that hybrid line-X orderings can work very well with asynchronous point-block ILU preconditioning.
\FloatBarrier

\section{Numerical results}
\label{sec:results}

We provide some experimental (numerical) results which demonstrate the performance of asynchronous block ILU preconditioning, albeit for simple cases. In what follows, `ILU' will always denote fixed-pattern incomplete LU factorization where the non-zero pattern of the $L$ and $U$ factors is the same as the original matrix; that is, no fill-in is allowed. Since our Jacobian matrices are naturally made up of small fixed-size dense blocks, the number of non-zero entries in the scalar and block representations are the same. Thus, the number of non-zero entries in the ILU(0) factorization and the block ILU(0) factorization are the same.

We are interested in investigating the improvement in convergence and/or parallel scaling brought about by the use of the asynchronous block ILU preconditioner when compared to the Chow-Patel preconditioner for two test cases, one with inviscid flow and the other with viscous laminar flow.
Secondly, we investigate the impact of the grid ordering on the convergence and parallel scaling of the different variants of the asynchronous ILU preconditioner.
Please note that for the nonlinear solver tests in sections \ref{sec:results_cyl} and \ref{sec:results_naca}, we averaged the data over three runs. We observed that the deviation is small (less than 3.3\% in all cases). 

The following parameters are used for all the runs.
\begin{itemize}
    \item Convergence tolerance for the non-linear problem is a relative drop of $10^{-6}$ in the energy residual. The maximum CFL number (with regard to equation \eqref{eqn:cfl}) is 10,000.
    \item The FGMRES \cite{gmresflexible} solver from the Portable, Extensible Toolkit for Scientific Computing (PETSc) \cite{petsc-user-ref} is used. The restart length is set to 30.
    At every non-linear iteration (pseudo-time step), linear solver convergence criterion is a residual ($\lVert \bld{A}\bld{x}-\bld{b}\rVert_2$) reduction by one order of magnitude, but only up to a maximum of 60 FGMRES iterations.
	\item For every non-linear solve, a spatially first-order accurate solver is used to obtain an initial solution for the second-order solver. This initial first-order solver starts from free-stream conditions and is converged by 1 order of magnitude of the energy residual. It uses a sequential ILU preconditioner so that all evaluation runs start from exactly the same initial solution. This initial solver is not timed - all reported results are from the second-order solve \footnote{`First order' and `second order' here refer to only the right-hand side fluxes. The Jacobian matrix is always that of the first-order discretization. For a second-order solve, however, the Jacobian matrices are computed at second-order accurate solution vectors.}.
    \item The initial guess for the asynchronous ILU factorization is the Jacobian matrix itself, while the initial guess for the triangular solves is the zero vector.
	\item For the point-block solver, the \texttt{BAIJMKL} matrix storage format is used for the Jacobian. This is PETSc's interface to block sparse matrix storage in the Intel Math Kernel Library (MKL). However, we do not use MKL in the preconditioners.
	The Eigen matrix library \cite{eigen} is used for vectorized small dense matrix operations in the point-block preconditioner.
	\item The CPU is an Intel Xeon Phi 7230 processor with 64 cores. High bandwidth memory was used exclusively for all the runs. Only one thread was used per core. Each thread was bound to one hyper-thread context of a core.
	The OpenMP chunk size was 384 in case of block preconditioners, and 4 times that, 1536, for the scalar preconditioners (the block size is 4).
	\item The Intel C/C++ compiler version 17 was used with \texttt{-O3 -xmic-avx512} as optimization flags.
	\item In the graphs in this section, curves labelled as using one thread or one core correspond to regular sequential preconditioning. For example, asynchronous ILU(0) factorization when run with one thread is exactly traditional ILU(0) factorization. Note that this need not be true for all implementations of asynchronous iterations because of different ways that SIMD instructions can be used, but it is true for the one presented here.
\end{itemize}

Variants of asynchronous ILU preconditioners are evaluated and compared based on the criteria given below.
\begin{itemize}
\item Convergence of the non-linear solver with respect to cumulative number of linear solver iterations.
\item Asynchronous ILU residual after 1 build sweep, over all pseudo-time steps. The reason for using 1 build sweep is discussed below, where we study the effect of using different numbers of sweeps.
\item Diagonal dominance of lower and upper triangular factors as a function of pseudo-time steps (non-linear iterations). Since we use an iterative method to apply the triangular factors, their diagonal dominance is an indicator of the stability of the triangular solve.
\item Strong scaling based on wall-clock time taken by all preconditioning operations until convergence of the non-linear problem. A good scalability will show both the parallel efficiency of the asynchronous kernels and also the strength of the preconditioner to reduce the number of iterations required.
\end{itemize}

Chow and Patel suggested \cite{ilu:chowpatel} scaling the original matrix symmetrically (using the same row and column scaling factors), resulting in unit diagonal entries, and factorizing this scaled matrix using asynchronous ILU iterations.
This was attempted for the two test cases considered in this article, but no significant difference was observed as a result. Effectiveness of preconditioner parameters in terms stability and speed was insensitive to such scaling. 
This was true irrespective of whether the block variant was used or not, and the grid ordering (among RCM, 1WD, line and line-hybrid orderings). 
The fact that a non-dimensional form of the Navier-Stokes equations is used may play a part in this observation.
The results presented in this work use the un-scaled matrix to compute the factorization, unless stated otherwise.

\subsubsection*{Number of asynchronous sweeps within the preconditioner}
\label{subsec:paramsweeps}

We now carry out a study to see how the general trend of FGMRES convergence of a linear system depends on the number of sweeps used to build and apply the asynchronous block ILU preconditioner. For this purpose, we extract the Jacobian matrix from an intermediate pseudo-time step for each test case and solve the linear system using the asynchronous block-ILU preconditioner and FGMRES(30) for two different thread settings - 16 and 62 threads.
The result for each sweep setting was averaged over 10 repeated runs for these tables. The maximum relative deviation over all sweep settings is reported in the caption of each table. Relative deviation is defined here for each sweep setting as the standard deviation divided by the average number of iterations over the 10 repetitions for that sweep setting.
Table \ref{tab:paramsweep-inv} shows the iteration counts for the inviscid cylinder case with reverse Cuthill-McKee (RCM) ordering, while table \ref{tab:paramsweep-visc-rcm} shows the same for the viscous NACA0012 airfoil case with RCM ordering and table \ref{tab:paramsweep-visc-1wd} corresponds to the viscous NACA0012 case with one-way dissection (1WD) ordering.

The first broad observation is that higher build sweeps typically require higher application sweeps for effective preconditioning. This is more pronounced for higher thread counts. However this trend can barely be seen for the 1WD ordering in case of viscous problems. This is the first indication that orderings such as RCM that typically work well for synchronous ILU may not work well for asynchronous ILU factorization.

\begin{table}[!htb]
\centering
\begin{tabular}{|c|c|c|c|c|c|c|c|c|c|c|c|c|c|c|}
    \hline
     & \multicolumn{7}{c|}{16 threads} & \multicolumn{7}{c|}{62 threads} \\
    \hline
    Apply sweeps & 1 & 2 & 3 & 5 & 10 & 20 & Exact & 1 & 2 & 3 & 5 & 10 & 20 & Exact \\
    Build sweeps & & & & & & & & & & & & & & \\
    \hline
    1 &      438  &  227  &  174  &  141  &  130  &  129  &  130  &  495  &  252  &  186  &  149  &  130   & 129  &  129 \\ 
    2 &      453  &  228  &  176  &  146  &  129  &  129  &  129  &  514  &  254  &  189  &  154  &  130   & 129  &  129 \\ 
    3 &      459  &  230  &  177  &  147  &  133  &  131  &  131  &  516  &  257  &  190  &  155  &  133   & 131  &  131 \\ 
    5 &      467  &  237  &  178  &  148  &  134  &  132  &  132  &  522  &  267  &  192  &  156  &  135   & 132  &  132 \\ 
   10 &      467  &  235  &  179  &  149  &  134  &  132  &  132  &  523  &  268  &  192  &  156  &  135   & 132  &  132 \\ 
   20 &      476  &  237  &  178  &  148  &  134  &  132  &  132  &  527  &  266  &  192  &  157  &  135   & 132  &  132 \\ 
Exact &      470  &  237  &  178  &  148  &  134  &  132  &  132  &  526  &  268  &  191  &  157  &  135   & 132  &  132 \\
    \hline
\end{tabular}
\caption{Number of FGMRES(30) iterations required for convergence to a relative tolerance of $1\times 10^{-2}$ as a function of number of sweeps used to build the asynchronous block-ILU preconditioner, for a matrix from the problem of inviscid flow over a cylinder (RCM ordering). Maximum deviation is about 2.2\%.}
\label{tab:paramsweep-inv}
\end{table}

\begin{table}[!htb]
\centering
\begin{tabular}{|c|c|c|c|c|c|c|c|c|c|c|c|c|c|c|}
    \hline
     & \multicolumn{7}{c|}{16 threads} & \multicolumn{7}{c|}{62 threads} \\
    \hline
    Apply sweeps & 1 & 2 & 3 & 5 & 10 & 20 & Exact & 1 & 2 & 3 & 5 & 10 & 20 & Exact \\
    Build sweeps & & & & & & & & & & & & & & \\
    \hline
     1   & 1525 & 471 &      &  109  &  107  &   61  &   45 &  &  &  &     & 334  & 99  & 53 \\
     2   &      &     &  134 &   30  &   41  &   20  &   22 &  &  &  & 150 &  35  & 30  & 21 \\
     3   &      &     & 1051 &   54  &   33  &   19  &   18 &  &  &  & 802 &  51  & 24  & 19 \\
     5   &  &  &  & 495 &  23 &   17  &   16 &  &  &  &  &  23 &   22  & 18 \\
    10   &  &  &  &     & 146 &   17  &   16 &  &  &  &  & 509 &   17  & 17 \\
    20   &  &  &  &     &     &   38  &   16 &  &  &  &  &     &   59  & 16 \\
  Exact  &  &  &  &     &     &  258  &   16 &  &  &  &  &     & 1147  & 16 \\
    \hline
\end{tabular}
\caption{Number of FGMRES(30) iterations required for convergence to a relative tolerance of $1\times 10^{-2}$ as a function of number of sweeps used to build the asynchronous block-ILU preconditioner, for a matrix from the problem of viscous flow over a NACA0012 airfoil. RCM ordering. Blanks indicate that the solver did not converge in 2500 iterations. Maximum relative deviation is 160\% in case of 1 build and 10 apply sweeps.}
\label{tab:paramsweep-visc-rcm}
\end{table}

\begin{table}[!htb]
\centering
\begin{tabular}{|c|c|c|c|c|c|c|c|c|c|c|c|c|c|c|}
    \hline
     & \multicolumn{7}{c|}{16 threads} & \multicolumn{7}{c|}{62 threads} \\
    \hline
    Apply sweeps & 1 & 2 & 3 & 5 & 10 & 20 & Exact & 1 & 2 & 3 & 5 & 10 & 20 & Exact \\
    Build sweeps & & & & & & & & & & & & & & \\
    \hline
        1 &  319  &  225  &  224  &  224  &  224  &  224  &  224  &  320  &  225  &  225  &  224  &  224  &  224  &  224 \\
        2 &  282  &  204  &  203  &  203  &  203  &  203  &  203  &  283  &  204  &  203  &  203  &  203  &  203  &  203 \\
        3 &  282  &  204  &  203  &  203  &  203  &  203  &  203  &  283  &  204  &  203  &  203  &  203  &  203  &  203 \\
        5 &  282  &  204  &  203  &  203  &  203  &  203  &  203  &  283  &  204  &  203  &  203  &  203  &  203  &  203 \\
       10 &  282  &  204  &  203  &  203  &  203  &  203  &  203  &  283  &  204  &  203  &  203  &  203  &  203  &  203 \\
       20 &  282  &  204  &  203  &  203  &  203  &  203  &  203  &  283  &  204  &  203  &  203  &  203  &  203  &  203 \\
    Exact &  282  &  204  &  203  &  203  &  203  &  203  &  203  &  283  &  204  &  203  &  203  &  203  &  203  &  203 \\
    \hline
\end{tabular}
\caption{Number of FGMRES(30) iterations required for convergence to a relative tolerance of $1\times 10^{-2}$ as a function of number of sweeps used to build the asynchronous block-ILU preconditioner, for a matrix from the problem of viscous flow over a NACA0012 airfoil. 1WD ordering. Maximum relative deviation is less than 1\%.}
\label{tab:paramsweep-visc-1wd}
\end{table}

For the inviscid case with RCM ordering with 62 threads, convergence generally gets worse if we use more build sweeps for a given number of application sweeps, while it improves if we use more application sweeps.
For the viscous case, RCM ordering results in poor performance of the asynchronous block-ILU preconditioner, and sensitivity to the number of sweeps is erratic (table \ref{tab:paramsweep-visc-rcm}).
It can be seen that the general trend is towards worse preconditioning with increasing build sweeps for constant apply sweeps, and better preconditioning for increasing apply sweeps with constant build sweeps. However, there are clearly exceptions - the erratic nature of this table supports the results shown later that RCM ordering leads to poor and unreliable performance for viscous cases.
Finally, the 1WD ordering gives robust results for this linear system of the viscous flow case. The results are almost independent of the number of threads, and all sweeps settings converge. For such cases, this ordering is clearly preferable to RCM, though the sequential (exact) preconditioning effectiveness is significantly worse. We further explore this in the context of convergence of the nonlinear problem further below (figure \ref{fig:abilu-orderings-resconv}).

From the tables, we see that there is no consistent and significant advantage of using more than 1 sweep to build the factorization. 
For the viscous case with 1WD ordering (table \ref{tab:paramsweep-visc-1wd}), there is an advantage to using two build sweeps, though it performs consistently with one build sweep as well.
Our objective is to use as few sweeps as possible without adversely impacting the convergence of the non-linear problem much.
Ultimately, for uniformity in the analysis of the non-linear solves in the next subsections, we choose 1 build sweep and 3 apply sweeps for all the studies.



\subsection{Inviscid subsonic flow over cylinder}
\label{sec:results_cyl}

This section will demonstrate the effectiveness of the asynchronous block ILU preconditioner for an inviscid flow over a cylinder. The unstructured mesh consists of 217,330 quadrilaterals. The free-stream Mach number is 0.38. The mesh is largely isotropic, therefore we do not use any line-based orderings for this case.

\begin{figure}[!htb]
	\begin{subfigure}{0.50\linewidth}
		\centering
	    \includegraphics[scale=0.52]{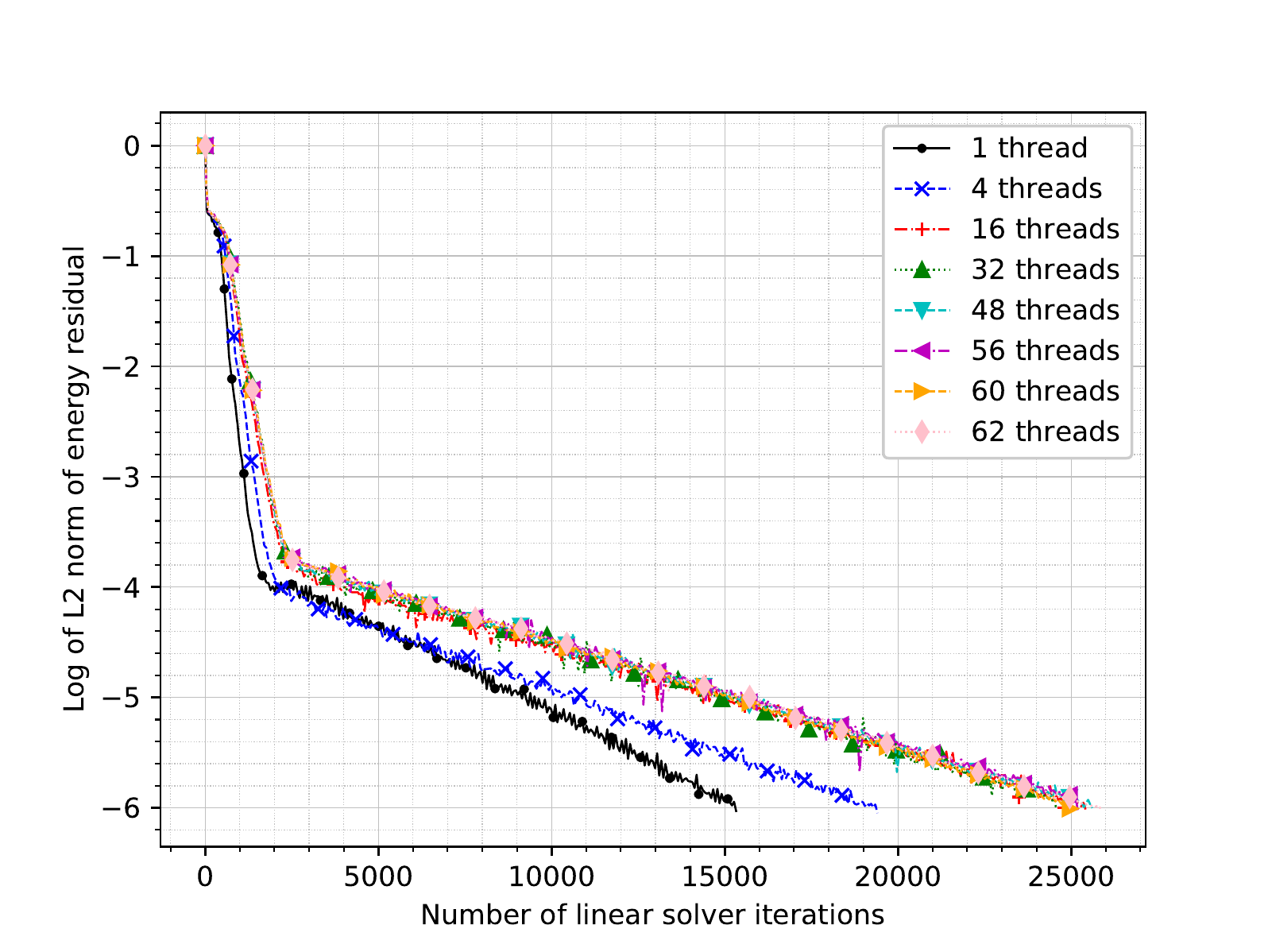}
		\subcaption{Async. scalar ILU, RCM ordering}
        \label{fig:inv-2dcyl-ailu-noscale-rcm-resconv}
	\end{subfigure}
	\begin{subfigure}{0.5\linewidth}
	\centering
	    \includegraphics[scale=0.52]{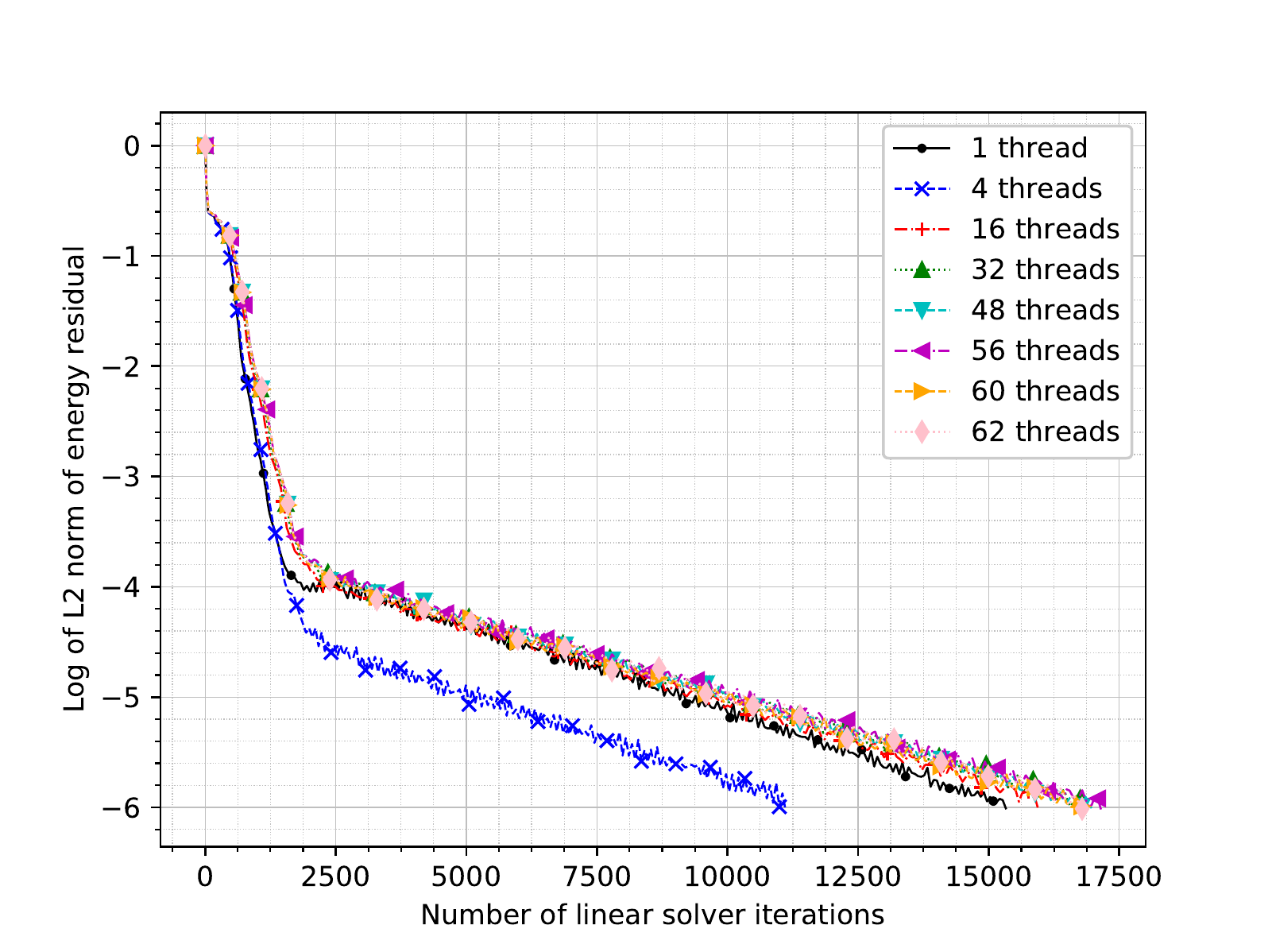}
	\subcaption{Async. block ILU, RCM ordering}
    \label{fig:inv-2dcyl-abilu-noscale-rcm-resconv}
	\end{subfigure}
	\begin{subfigure}{0.50\linewidth}
		\centering
	    \includegraphics[scale=0.52]{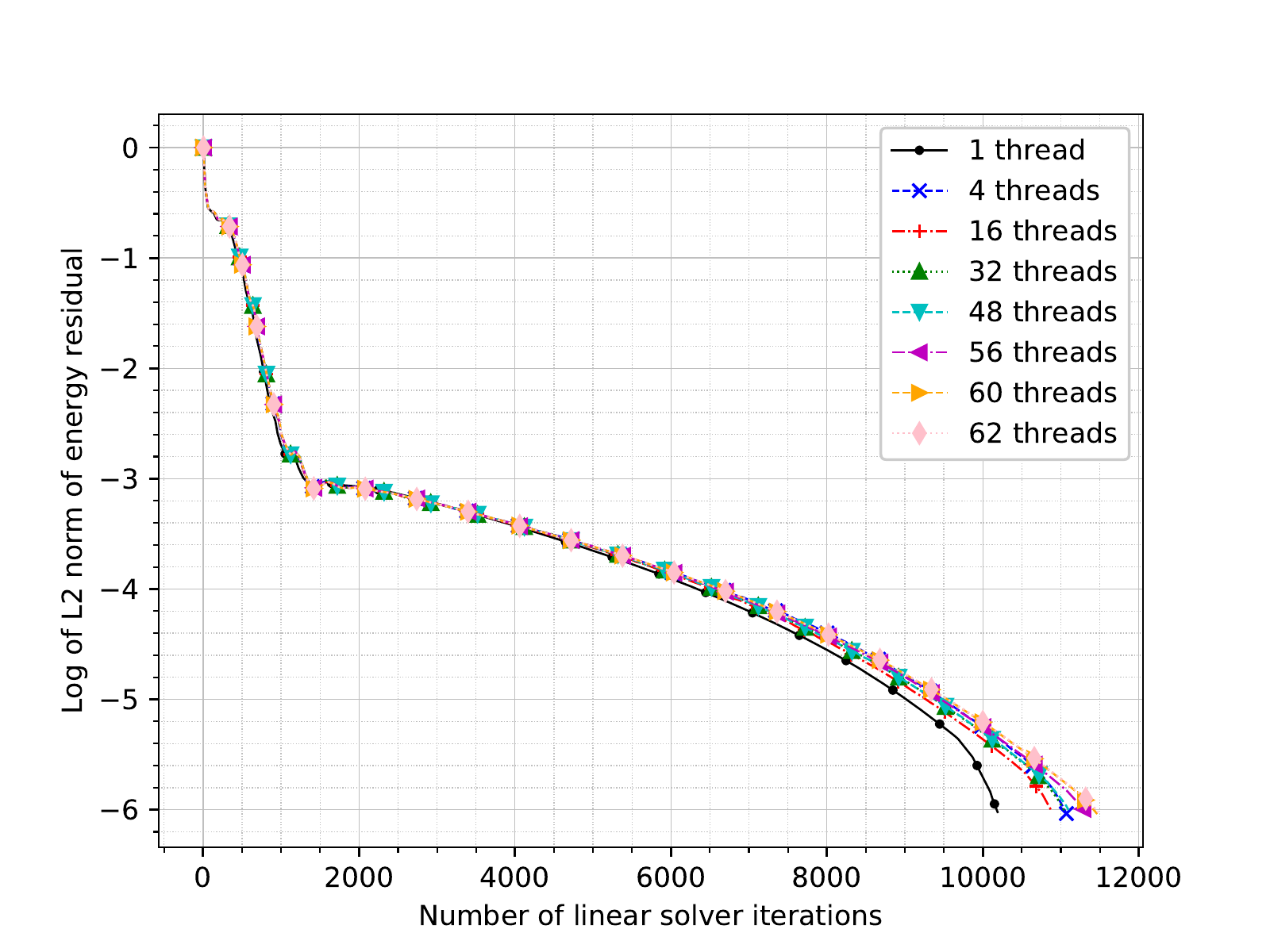}
		\subcaption{Async. scalar ILU, 1WD ordering}
    \label{fig:inv-2dcyl-ailu-noscale-1wd-resconv}
	\end{subfigure}
	\begin{subfigure}{0.5\linewidth}
	\centering
	    \includegraphics[scale=0.52]{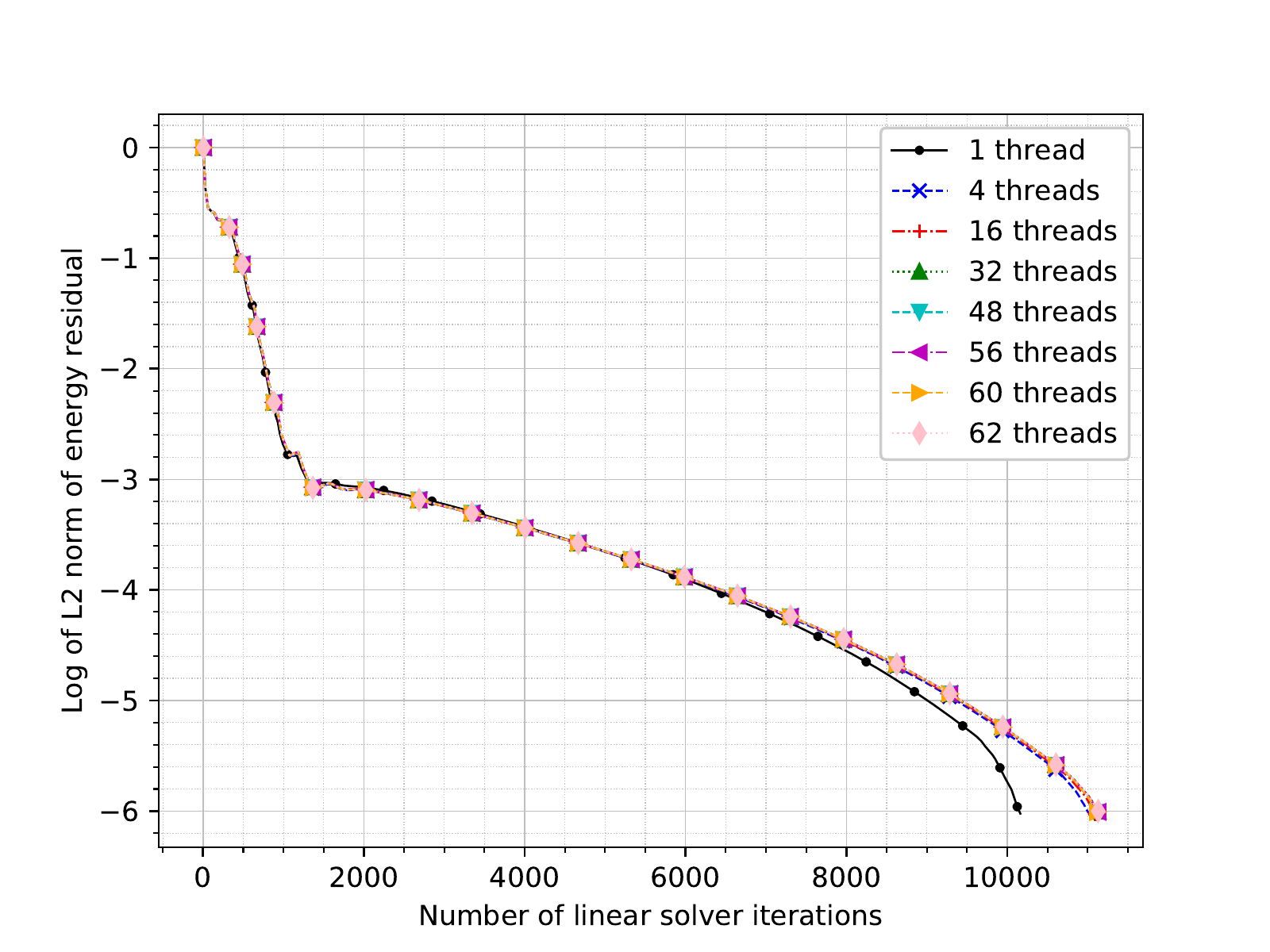}
	\subcaption{Async. block ILU, 1WD ordering}
    \label{fig:inv-2dcyl-abilu-noscale-1wd-resconv}
	\end{subfigure}
    \caption{Convergence of the non-linear problem w.r.t. cumulative FGMRES iterations, for the inviscid cylinder case, showing the advantage of asynchronous block ILU over the scalar variant}
    \label{fig:inv-2dcyl-blocking_ordering-resconv}
\end{figure}

\begin{figure}[!htb]
\begin{subfigure}{0.5\linewidth}
    \centering
    \includegraphics[scale=0.6]{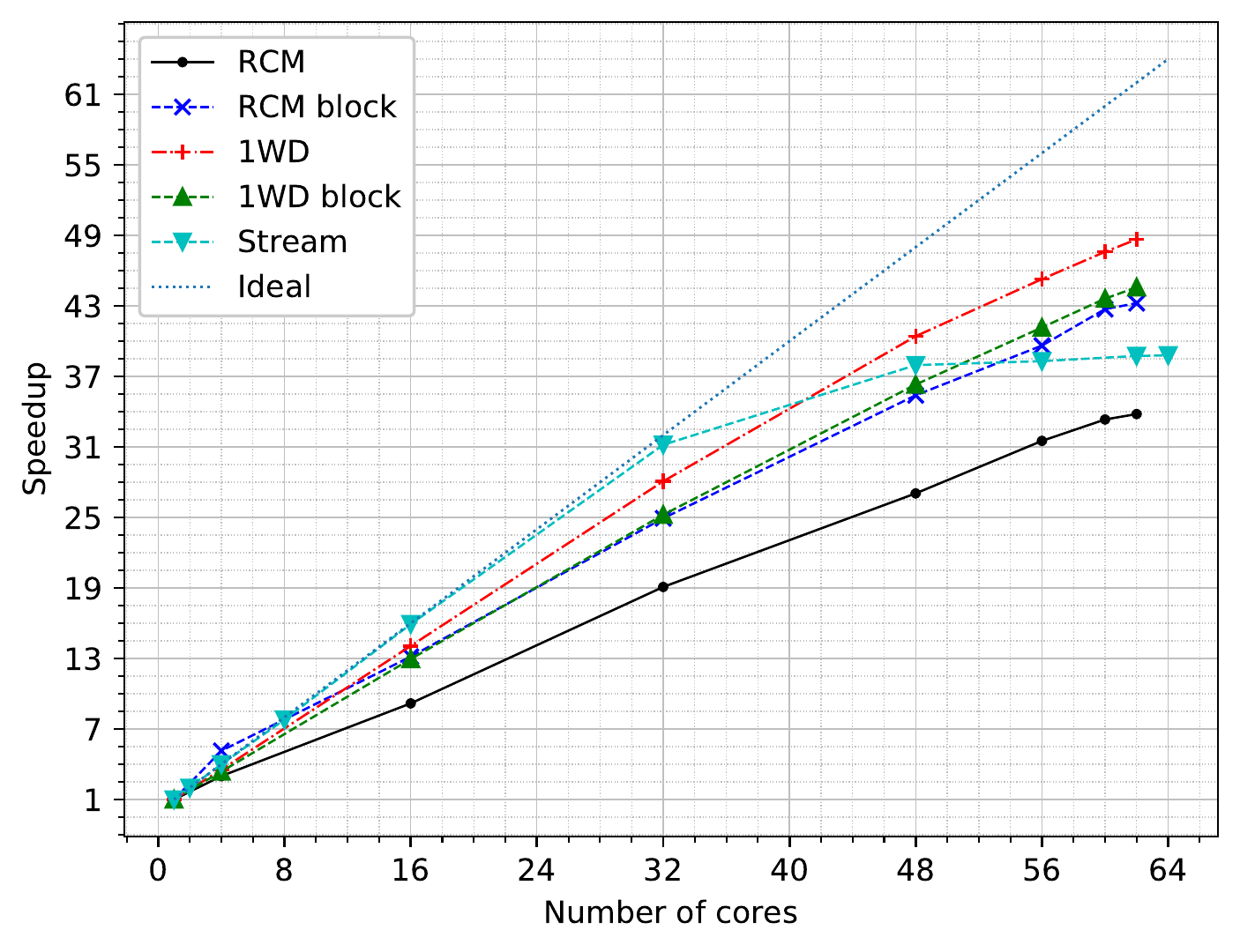}
    \caption{Strong scaling}
    \label{fig:inv-2dcyl-strongscaling}
\end{subfigure}
\begin{subfigure}{0.5\linewidth}
    \centering
    \includegraphics[scale=0.6]{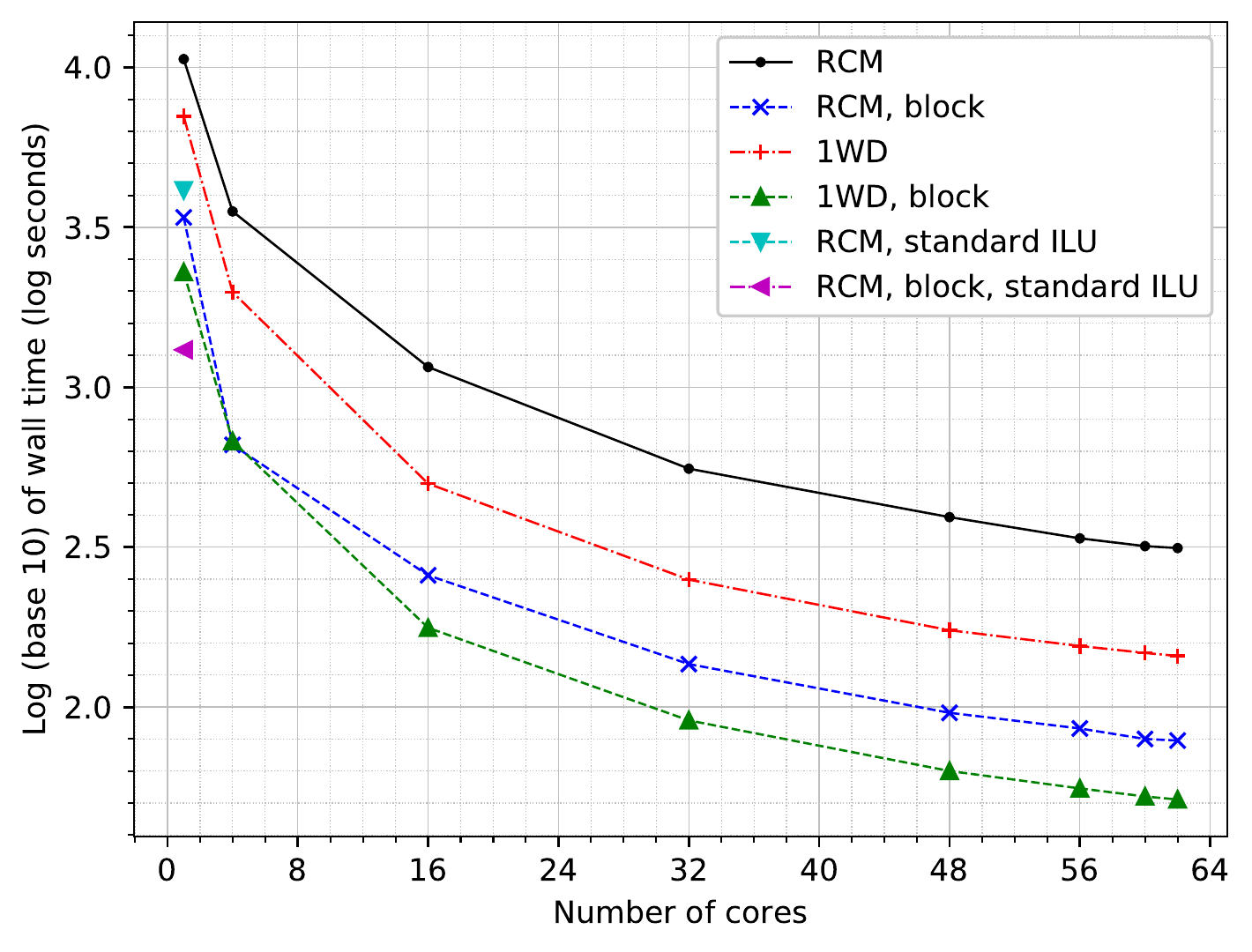}
    \caption{Wall-clock time}
    \label{fig:inv-2dcyl-walltimes}
\end{subfigure}
\caption{Strong scaling and wall-clock time taken by asynchronous scalar and block ILU preconditioning operations over the entire non-linear solve, using RCM and 1WD orderings, for the inviscid flow case}
\label{fig:inv-2dcyl-performance}
\end{figure}

First, we look at the convergence of the entire non-linear CFD problem, in terms of the norm of the energy residual, with respect to the cumulative number of FGMRES iterations in figure \ref{fig:inv-2dcyl-blocking_ordering-resconv}. 
In case of RCM ordering of the grid cells (figure \ref{fig:inv-2dcyl-ailu-noscale-rcm-resconv}), we observe that the asynchronous scalar ILU preconditioner causes a significant increase in the required number of FGMRES iterations as we increase the number of threads. Asynchronous block ILU, however, leads to a much smaller increase as we increase the number of threads (figure \ref{fig:inv-2dcyl-abilu-noscale-rcm-resconv}). This is reflected in the strong scaling and wall-clock time plots shown later (figures \ref{fig:inv-2dcyl-performance}).
For this case, with the RCM ordering, the 4-thread run behaves anomalously for asynchronous block ILU - it converges even faster than the sequential ILU preconditioner. We have sometimes seen such abnormally fast convergence in previous work as well \cite{async:aditya_absgs_2019}, but we do not generally expect an asynchronous ILU preconditioner to converge faster (in terms of number of linear solver iterations) than the corresponding standard sequential preconditioner.

However, when applied after one-way dissection (1WD) ordering, asynchronous ILU shows no such effects. In terms of the required number of FGMRES iterations for convergence, there is no significant difference between the scalar- and block-ILU preconditioners, though the block preconditioner does appear to be slightly less sensitive to the number of threads in the final few iterations (compare figures \ref{fig:inv-2dcyl-ailu-noscale-1wd-resconv} and \ref{fig:inv-2dcyl-abilu-noscale-1wd-resconv}). We will see that the asynchronous block ILU preconditioner is faster in terms of wall-clock time (figure \ref{fig:inv-2dcyl-walltimes}). Additionally, looking at the number of FGMRES iterations along the $x$-axis, we note that this ordering gives a much better ILU preconditioner than the RCM ordering for this test case.

The difference between scalar and block ILU can also be described by the variation of the ILU fixed-point residual and the minimum diagonal dominance of the lower and upper triangular factors as the nonlinear solve proceeds.
We observed that one sweep of asynchronous block ILU leads to a relatively more accurate factorization compared to the original scalar variant. 
Furthermore, the block $L$-factor has a lower max-norm for its Jacobi iteration matrix, which means it is relatively more diagonally dominant in the block case than in the scalar case. The max-norm spikes to very high levels at some time steps for the scalar $L$-factor. The data for this has been omitted for the sake of brevity, but we return to this consideration in the next subsection for the viscous NACA0012 case.


Finally, figure \ref{fig:inv-2dcyl-performance} shows the performance of the asynchronous ILU variants in terms of wall-clock time. Note that each data point in these graphs represents the time taken by all preconditioning operations over the entire non-linear solver. Thus, any slowdowns because of a weaker preconditioner are represented here. We also include the scaling of the Stream benchmark \cite{stream} to compare against the scaling of a strongly memory bandwidth limited code.
In the speedup plot (figure \ref{fig:inv-2dcyl-strongscaling}), we observe that the scalar asynchronous ILU with RCM ordering does not show good parallel scaling. For all other variants, the parallel scaling is favourable compared to Stream. 
Since our kernels have more floating-point operations and less trivial memory access patterns than Stream, they are unlikely to be as memory-bandwidth limited. This is supported by the fact that some of the asynchronous ILU variants continue strong scaling after Stream reaches its limit.
In the wall-clock time plot (figure \ref{fig:inv-2dcyl-walltimes}), we see that asynchronous block ILU with 1WD ordering performs the best at high core counts. We include one data point each for the sequential ILU and block ILU preconditioners for reference.

\subsection{Viscous laminar flow over NACA0012 airfoil}
\label{sec:results_naca}

We now investigate a case of viscous laminar flow over a NACA0012 airfoil. The grid is two-dimensional unstructured, made up of quadrilaterals in the boundary layer and triangles elsewhere, obtained from the SU2 \cite{cfdsolvers:su2} test case repository \cite{su2testcases}. The Mach number is 0.5 and Reynolds number is 5000.
The total number of cells is 210,496, so the dimension of the problem is 841,984.

\begin{figure}[!htb]
	\begin{subfigure}{0.52\linewidth}
		\centering
	    \includegraphics[scale=0.52]{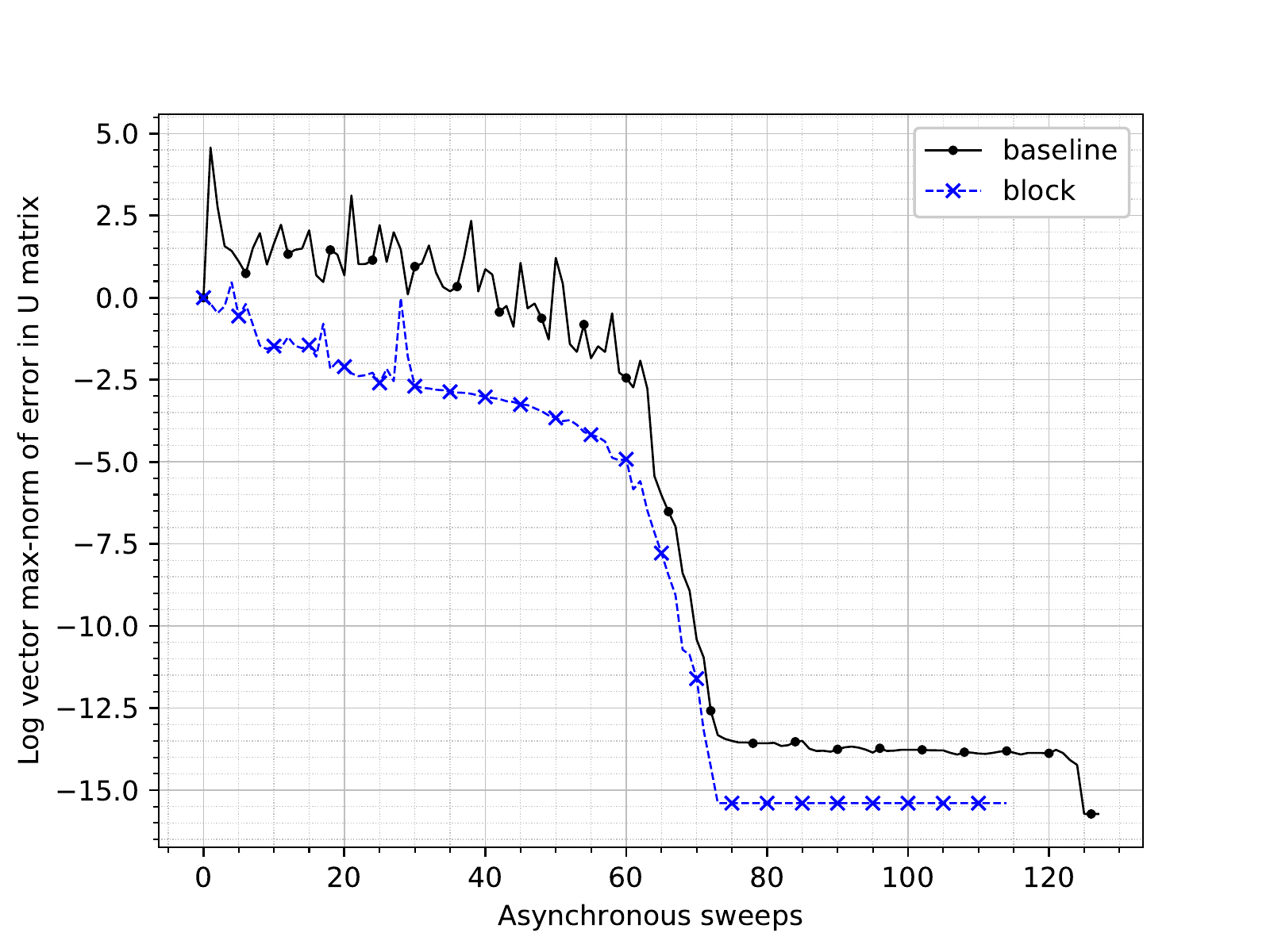}
		\subcaption{$\lVert \mathtt{vec}(U_{async}-U_{ilu})\rVert_\infty/\lVert \mathtt{vec}(U_0-U_{ilu})\rVert_\infty$}
    \label{fig:iluconv-rcm-upperdiff}
	\end{subfigure}
	\begin{subfigure}{0.52\linewidth}
	\centering
	    \includegraphics[scale=0.52]{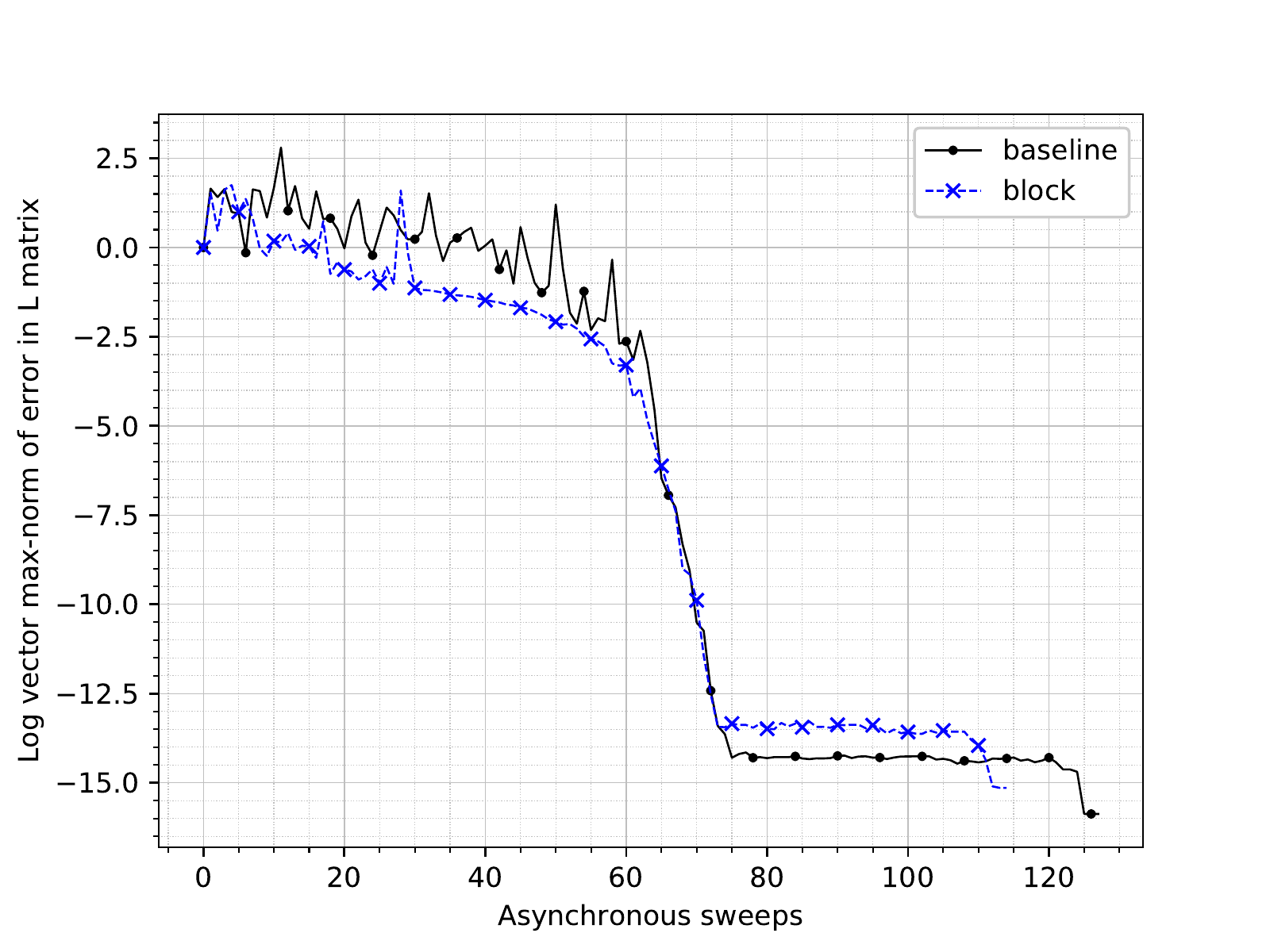}
		\subcaption{$\lVert \mathtt{vec}(L_{async}-L_{ilu})\rVert_\infty/\lVert \mathtt{vec}(L_0-L_{ilu})\rVert_\infty$}
    \label{fig:iluconv-rcm-lowerdiff}
	\end{subfigure}
    \caption{Convergence of different async. ILU fixed-point iterations (RCM ordering, 62 threads) for the viscous flow case (averaged over 5 runs)}
    \label{fig:iluconv-visc-naca0012}
\end{figure}

We first demonstrate the convergence of the asynchronous ILU sweeps to the `exact' incomplete (ILU(0)) $L$ and $U$ factors for one of the linear systems required for solving the non-linear problem.
The grid has been ordered in the RCM ordering and 62 threads have been used (figure \ref{fig:iluconv-visc-naca0012}). The matrix is taken from a pseudo-time step at which the CFL number is 4241 and the energy residual is $1.5\times 10^{-4}$.
The `baseline' run is a scalar asynchronous ILU(0) iteration. The errors in the $L$- and $U$-factors are normalized by the initial error.
We see that all the iterations converge to machine precision after a sufficient number of sweeps. We also see the asymptotically near-instantaneous convergence because the spectral radius of the iterations' Jacobian matrix is zero. Thus we regard this as numerical evidence not only of global convergence (theorem \ref{thm:blockilu_globalconv})  but also of asymptotically trivial local convergence (theorem \ref{thm:async_blockilu_localconv}).

\begin{figure}[!htb]
	\begin{subfigure}{0.52\linewidth}
		\centering
	    \includegraphics[scale=0.52]{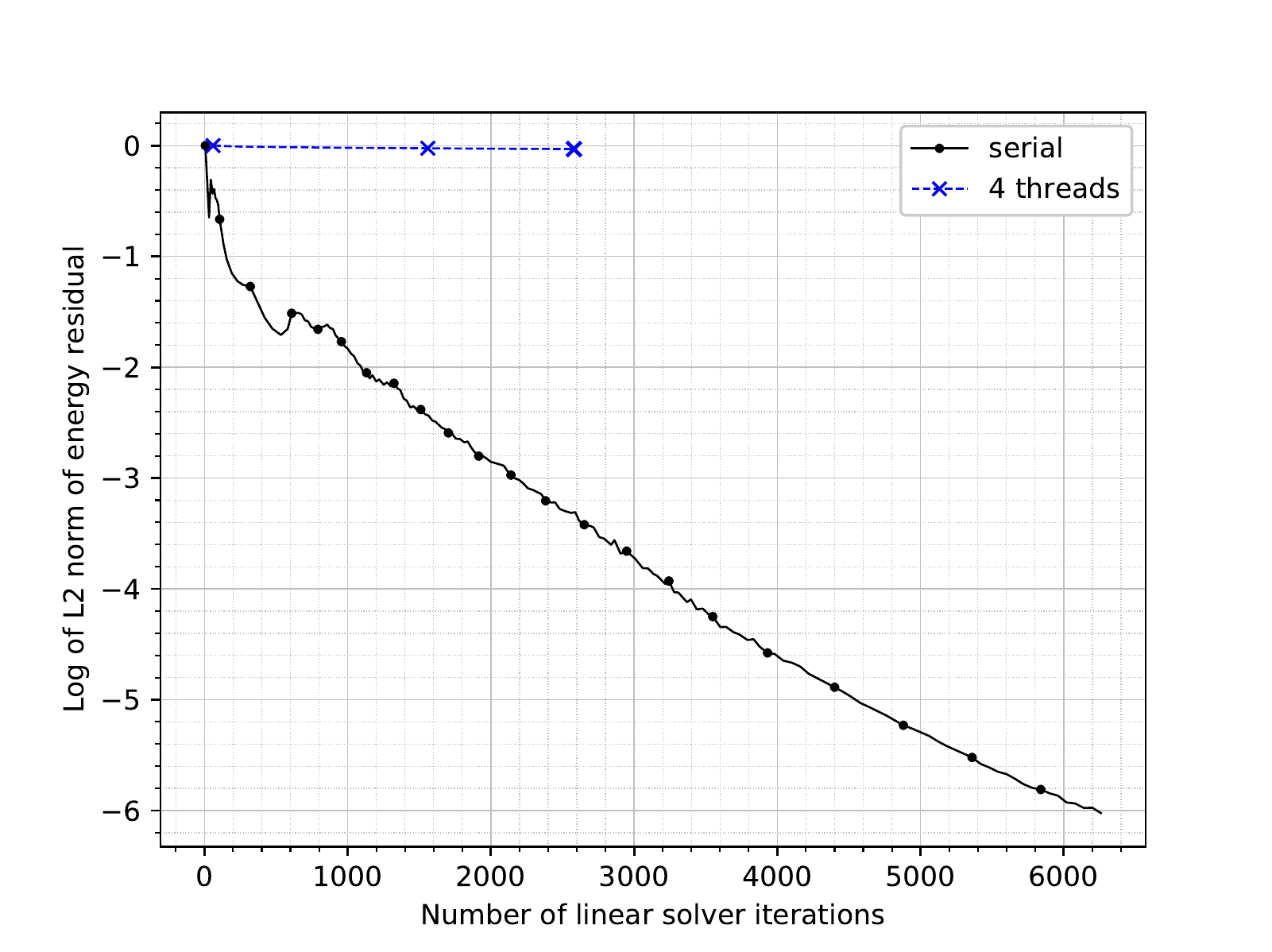}
		\subcaption{Async. scalar ILU}
    \label{fig:ailu-scaled-line_1wd-resconv}
	\end{subfigure}
	\begin{subfigure}{0.52\linewidth}
	\centering
	    \includegraphics[scale=0.52]{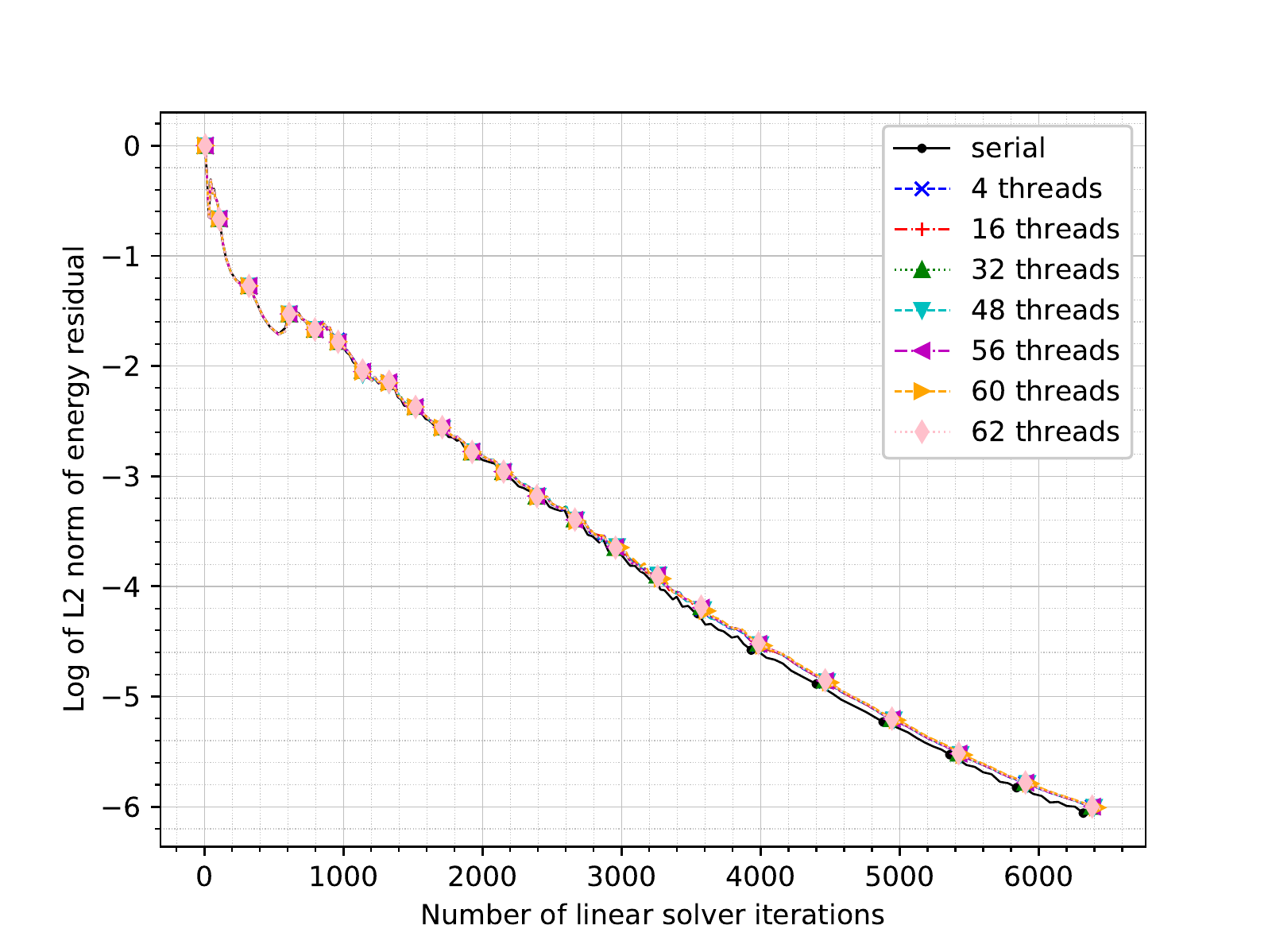}
	\subcaption{Async. block ILU}
    \label{fig:abilu-scaled-line_1wd-resconv}
	\end{subfigure}
    \caption{Convergence (or lack thereof) of the non-linear problem w.r.t. cumulative FGMRES iterations; line-1WD ordering and with symmetric scaling}
\end{figure}

Next, we look at the convergence of the non-linear problem, in terms of the norm of the energy residual, with respect to the cumulative number of FGMRES iterations.
The first thing to note is that the scalar asynchronous ILU(0) solver does not work for this case whenever more than 1 thread is used,
even when the matrix is first scaled symmetrically (figure \ref{fig:ailu-scaled-line_1wd-resconv}).
Thus, in this case, the block variant is necessary for obtaining convergence with asynchronous ILU(0) preconditioning (figure \ref{fig:abilu-scaled-line_1wd-resconv}).

We attempt to provide some insight into this difference in convergence by showing the variation of the ILU residual and the minimum diagonal dominance of the lower and upper triangular factors.
\begin{figure}[!htb]
	\begin{subfigure}{0.52\linewidth}
		\centering
	    \includegraphics[scale=0.52]{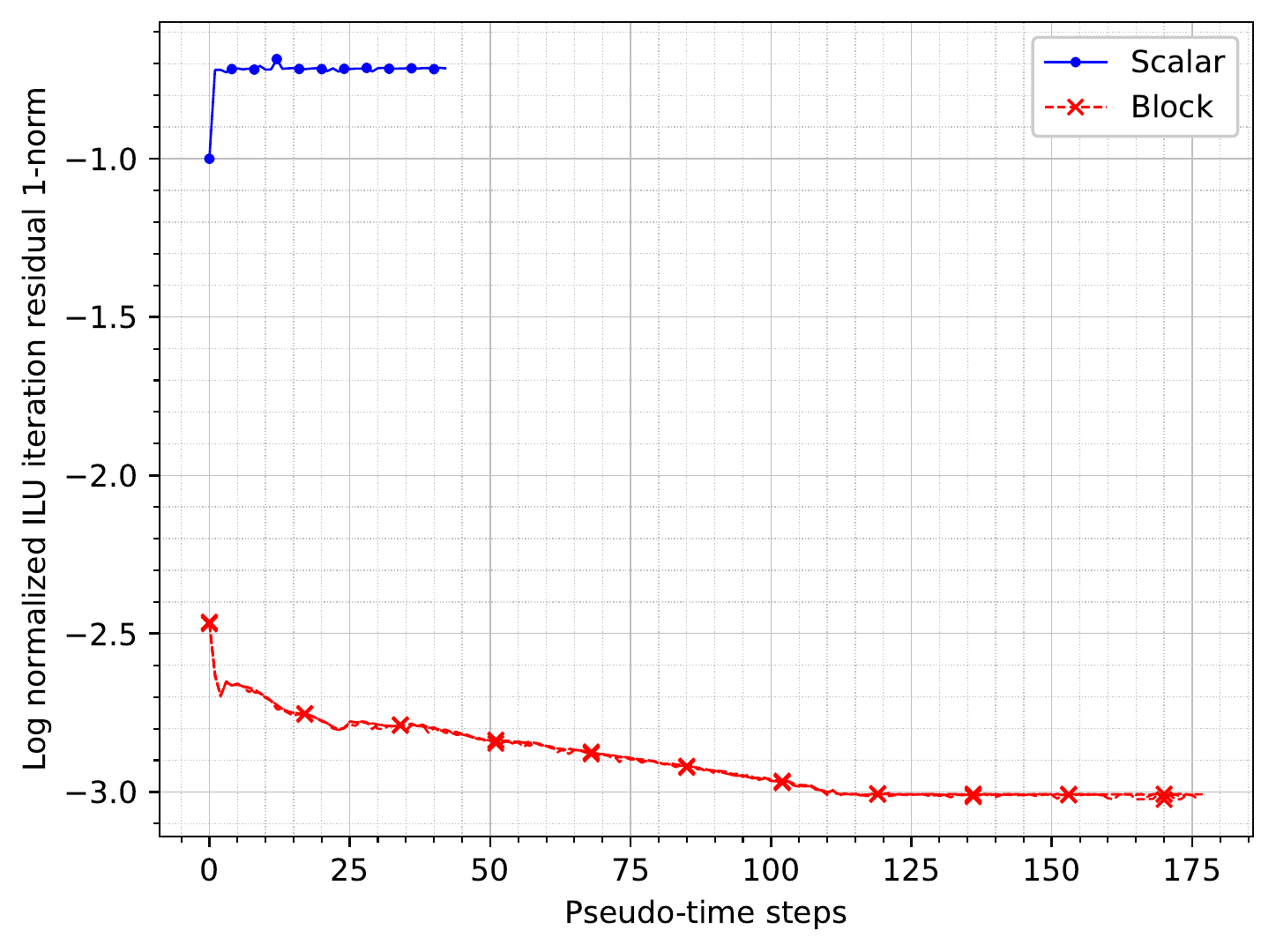}
		\subcaption{Reduction in ILU fixed-point residual vector norm}
        \label{fig:ilures-ailu-scaled-line_1wd}
	\end{subfigure}
	\begin{subfigure}{0.52\linewidth}
	\centering
	    \includegraphics[scale=0.52]{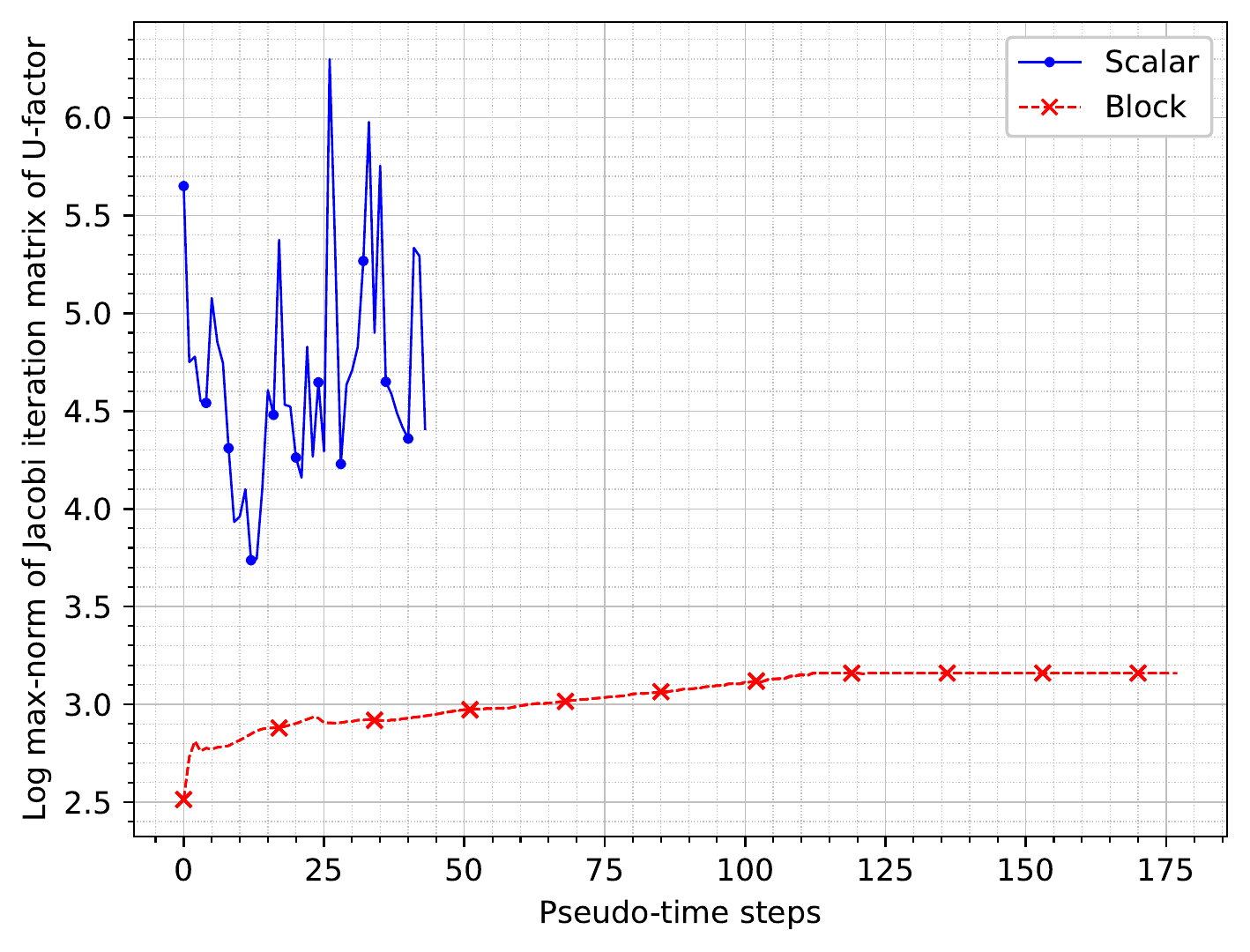}
	\subcaption{Norm of Jacobi iteration matrix for $U$-factor}
    \label{fig:udom-scaled-line_1wd}
	\end{subfigure}
	\begin{subfigure}{0.52\linewidth}
	\centering
	    \includegraphics[scale=0.52]{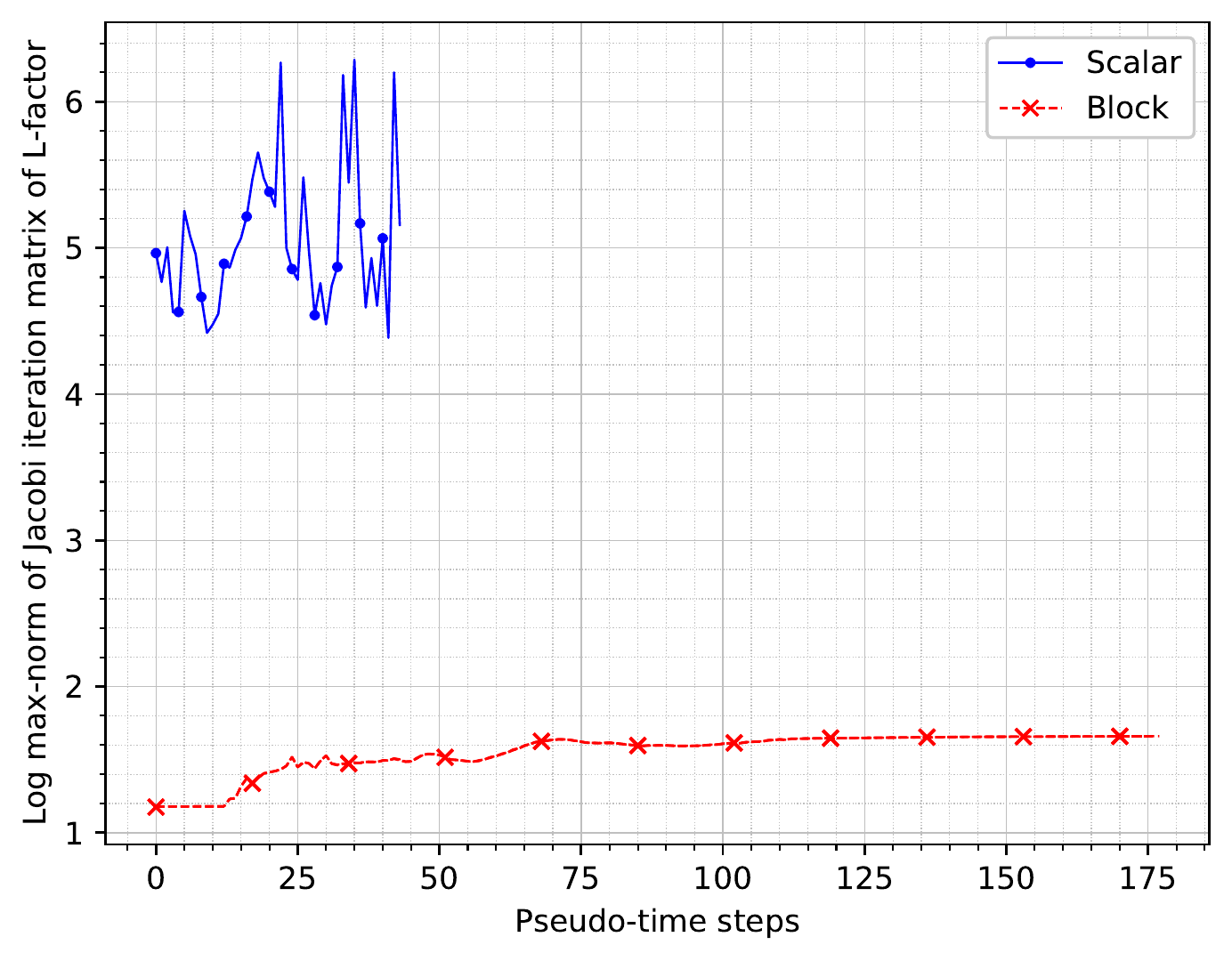}
	\subcaption{Norm of Jacobi iteration matrix for $L$-factor}
    \label{fig:ldom-ailu-scaled-line_1wd}
	\end{subfigure}
	\begin{subfigure}{0.52\linewidth}
    	\centering
	    \includegraphics[scale=0.52]{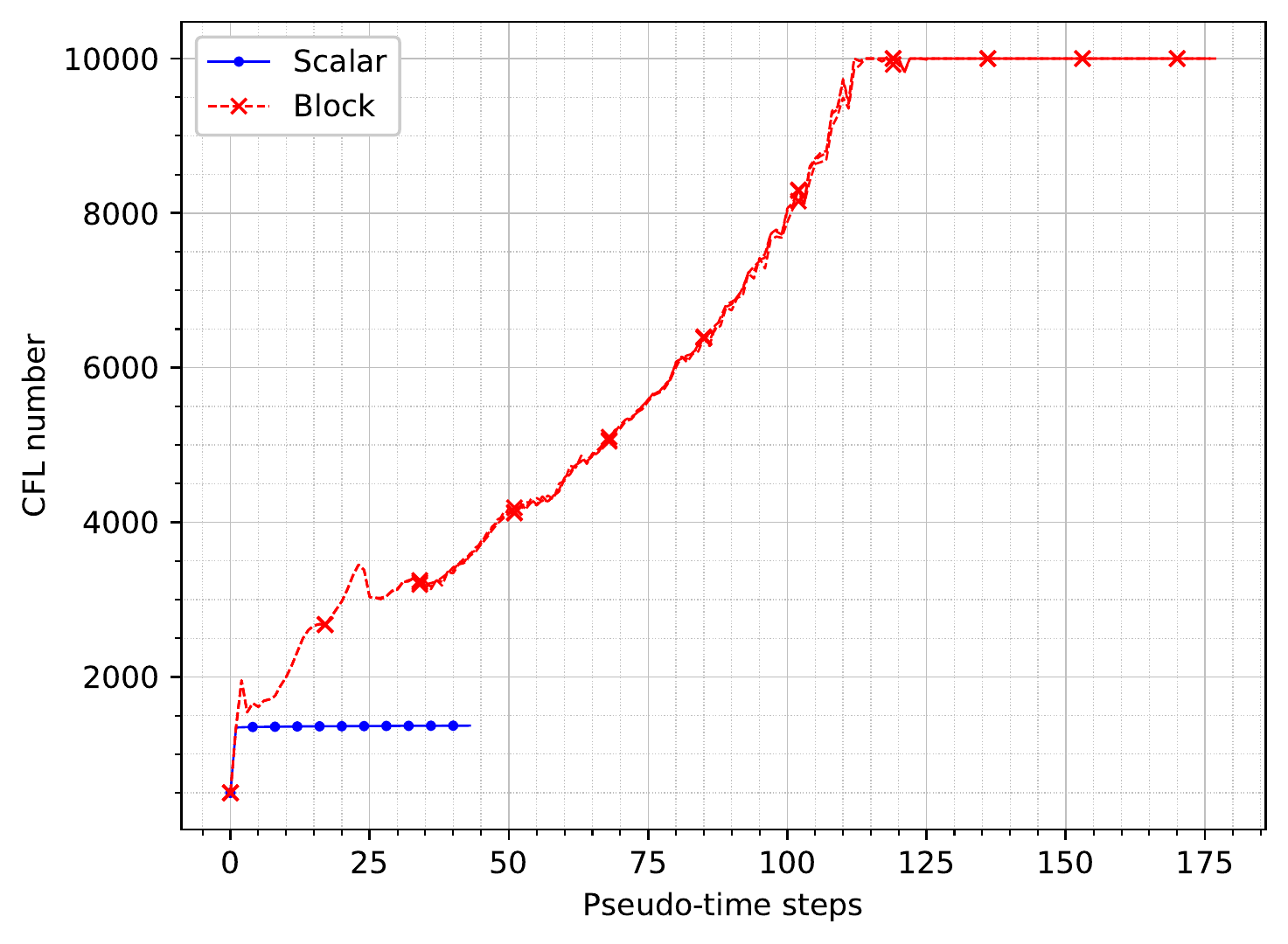}
    	\subcaption{CFL number}
        \label{fig:cfl-ailu-scaled-line_1wd}
	\end{subfigure}
    \caption{Properties of async. ILU factorization, and CFL number, w.r.t. pseudo-time steps for the viscous flow case with line-1WD ordering; 4 threads}
    \label{fig:scalar_vs_blk-precinfo}
\end{figure}
It can be seen that the ILU factorization is more accurate in the block case (in the vector 1-norm $\lVert \bld{x}_1 - \bld{g}_{(b)ilu}(\bld{x}_1) \rVert_1/\lVert \bld{x}_0 - \bld{g}_{(b)ilu}(\bld{x}_0) \rVert_1$, figure \ref{fig:ilures-ailu-scaled-line_1wd}).
Further, the computed $L$ and $U$ factors have better diagonal dominance property (figures \ref{fig:udom-scaled-line_1wd}, \ref{fig:ldom-ailu-scaled-line_1wd}), which is shown by the lower max-norm of the Jacobi iteration matrix for the lower and upper triangular solves. Lower max-norm of the Jacobi iteration matrix is equivalent to better diagonal dominance. (The $L$ and $U$ matrices are not actually diagonally dominant in either case - the Jacobi iteration matrix max-norm is greater than 1 in both cases - but the block case is better in this regard.) We also observe that the norm of the Jacobi iteration matrix for ILU factors increases moderately until the CFL number increases.
Recall that the CFL number is adjusted every non-linear iteration based on the ratio of the current and previous residual norms, as described in section \ref{sec:solver}.

\begin{figure}[!htb]
	\begin{subfigure}{0.5\linewidth}
		\centering
	    \includegraphics[scale=0.52]{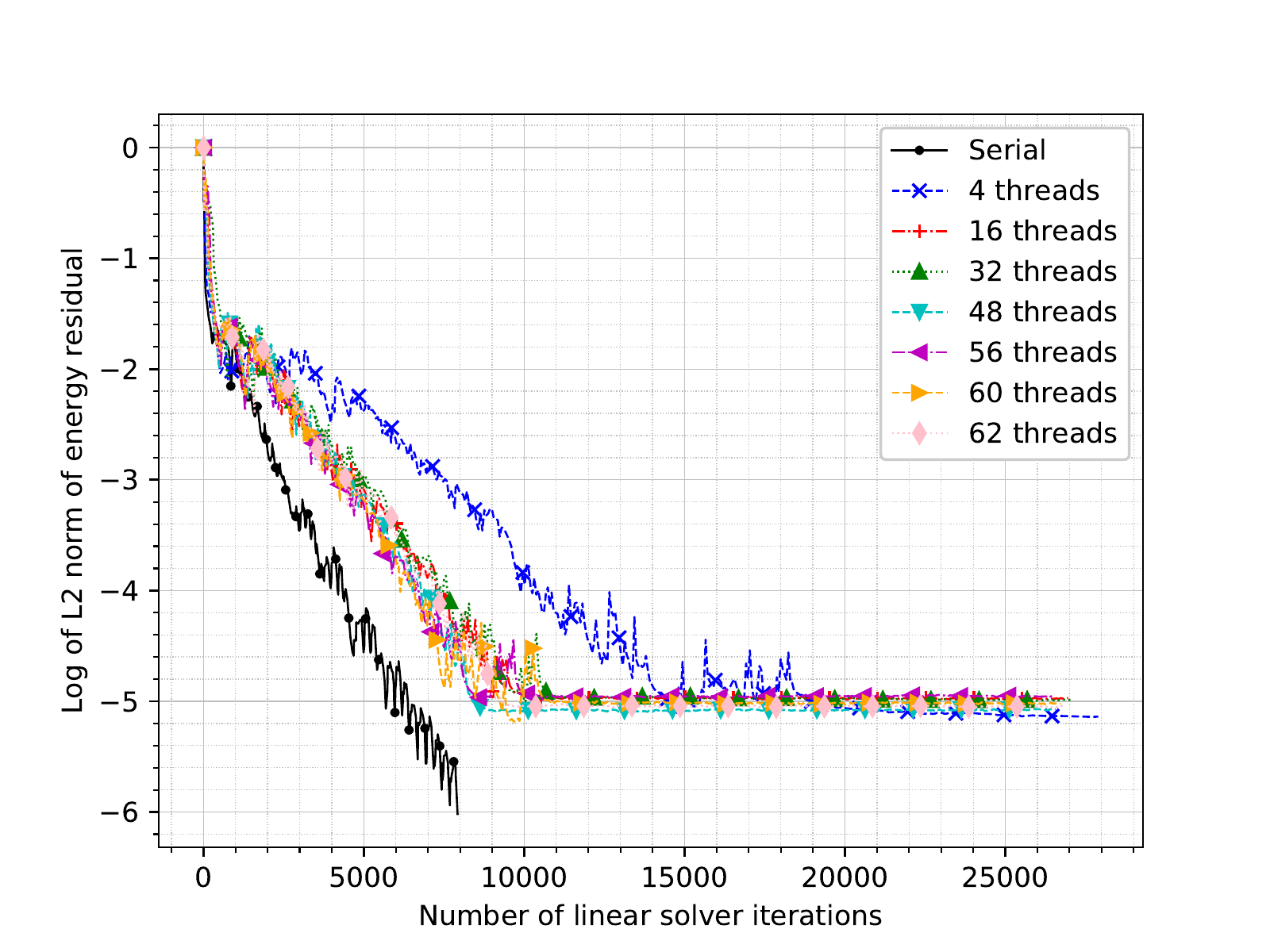}
	    \subcaption{RCM ordering}
    \label{fig:abilu-noscale-rcm-resconv}
	\end{subfigure}
	\begin{subfigure}{0.5\linewidth}
		\centering
	    \includegraphics[scale=0.52]{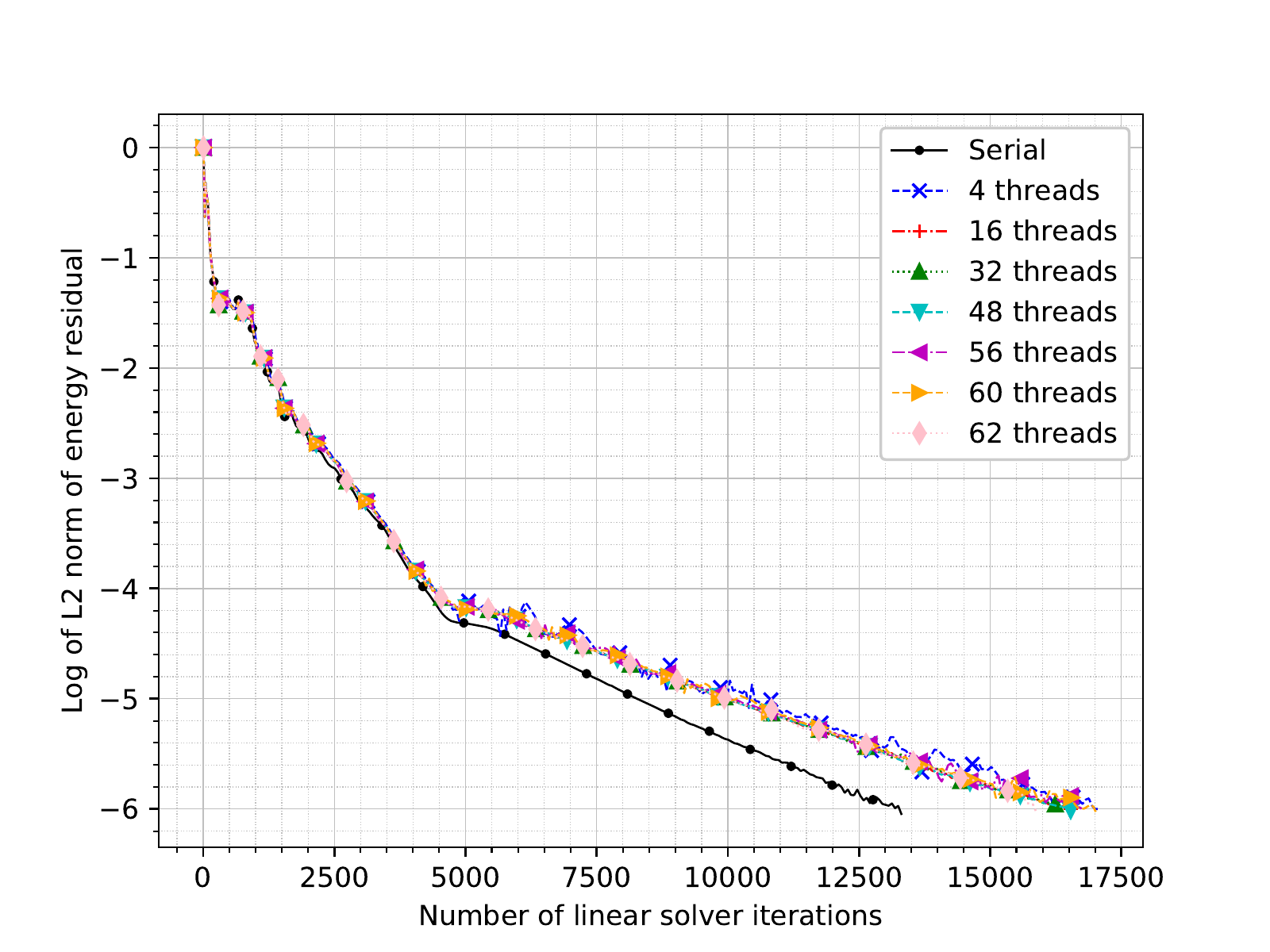}
	    \subcaption{1WD ordering}
    \label{fig:abilu-noscale-1wd-resconv}
	\end{subfigure}
	\begin{subfigure}{0.5\linewidth}
	\centering
	    \includegraphics[scale=0.52]{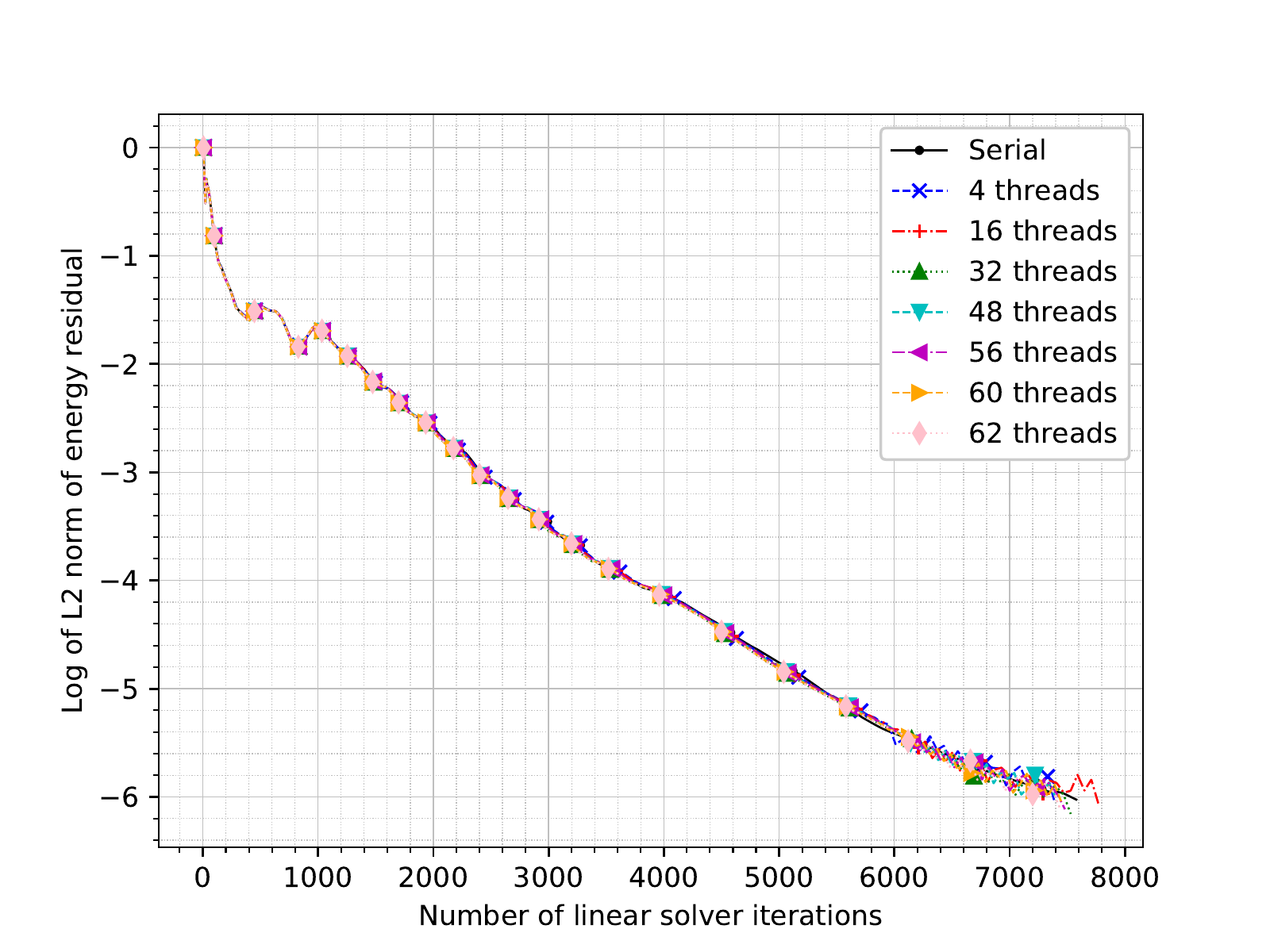}
		\subcaption{Line ordering}
    \label{fig:abilu-noscale-line-resconv}
	\end{subfigure}
	\begin{subfigure}{0.5\linewidth}
	\centering
	    \includegraphics[scale=0.52]{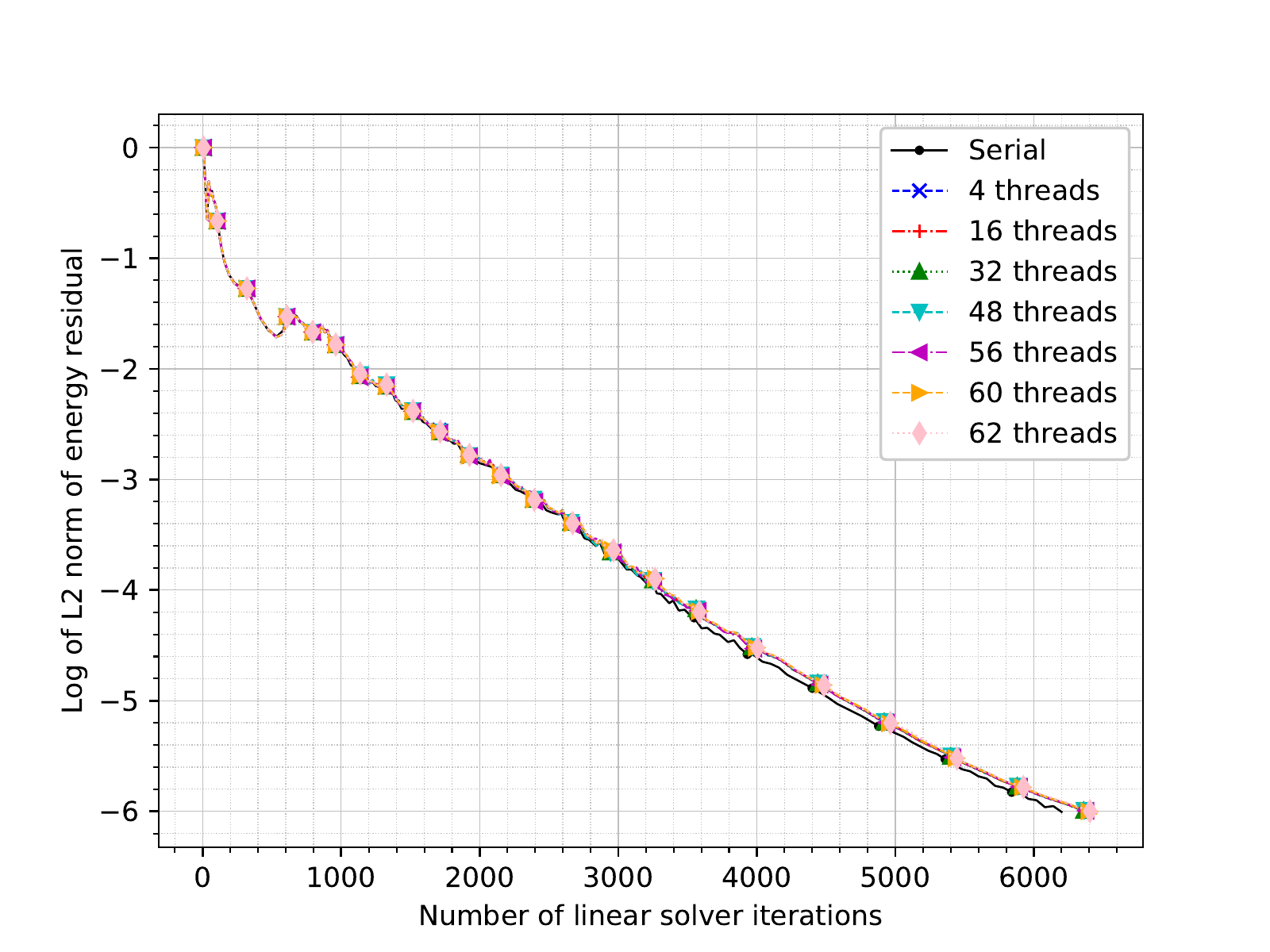}
		\subcaption{Line 1WD ordering}
    \label{fig:abilu-noscale-line_1wd-resconv}
	\end{subfigure}
    \caption{Convergence of non-linear problem w.r.t. cumulative FGMRES iterations with different orderings (without scaling of the original matrices)}
    \label{fig:abilu-orderings-resconv}
\end{figure}

Next, the asynchronous block ILU solver with RCM ordering stalls for any more than 1 thread (figure \ref{fig:abilu-noscale-rcm-resconv}).
The one-way dissection (1WD) ordering is able to recover convergence, but the parallel runs still require more FGMRES iterations than the serial run (figure \ref{fig:abilu-noscale-1wd-resconv}).
Finally, the line-based orderings converge in an essentially thread-independent number of iterations (figures \ref{fig:abilu-noscale-line-resconv}, \ref{fig:abilu-noscale-line_1wd-resconv}). The line-1WD hybrid ordering converges in the least number of iterations.

\begin{figure}[!htb]
	\begin{subfigure}{0.52\linewidth}
		\centering
	    \includegraphics[scale=0.52]{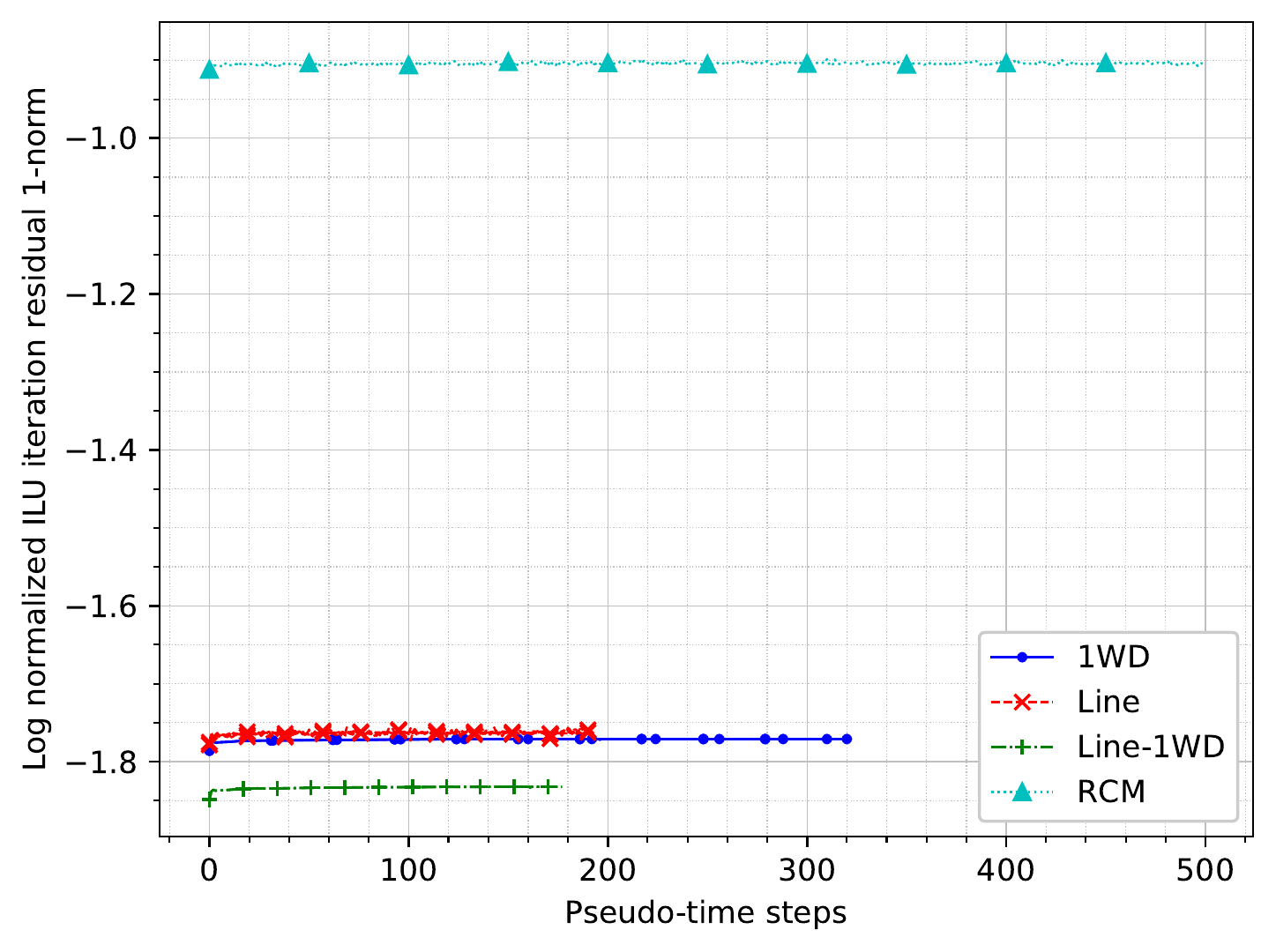}
		\subcaption{ILU residual vector norm}
        \label{fig:ilures-abilu-noscale-orderings}
	\end{subfigure}
	\begin{subfigure}{0.52\linewidth}
	\centering
	    \includegraphics[scale=0.52]{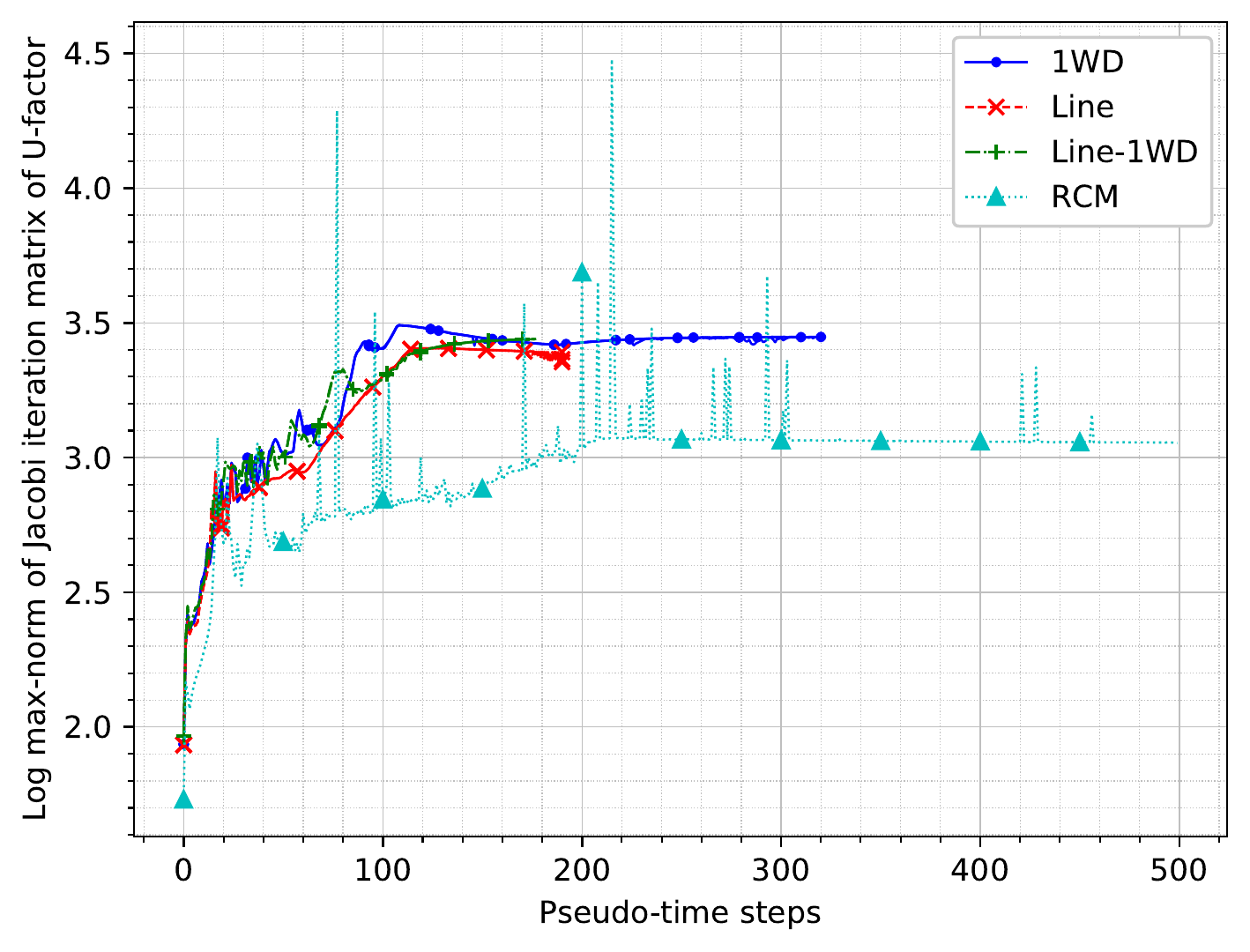}
	\subcaption{Norm of Jacobi iteration matrix for $U$-factor}
    \label{fig:udom-abilu-noscale-orderings}
	\end{subfigure}
	\begin{subfigure}{0.51\linewidth}
	\centering
	    \includegraphics[scale=0.52]{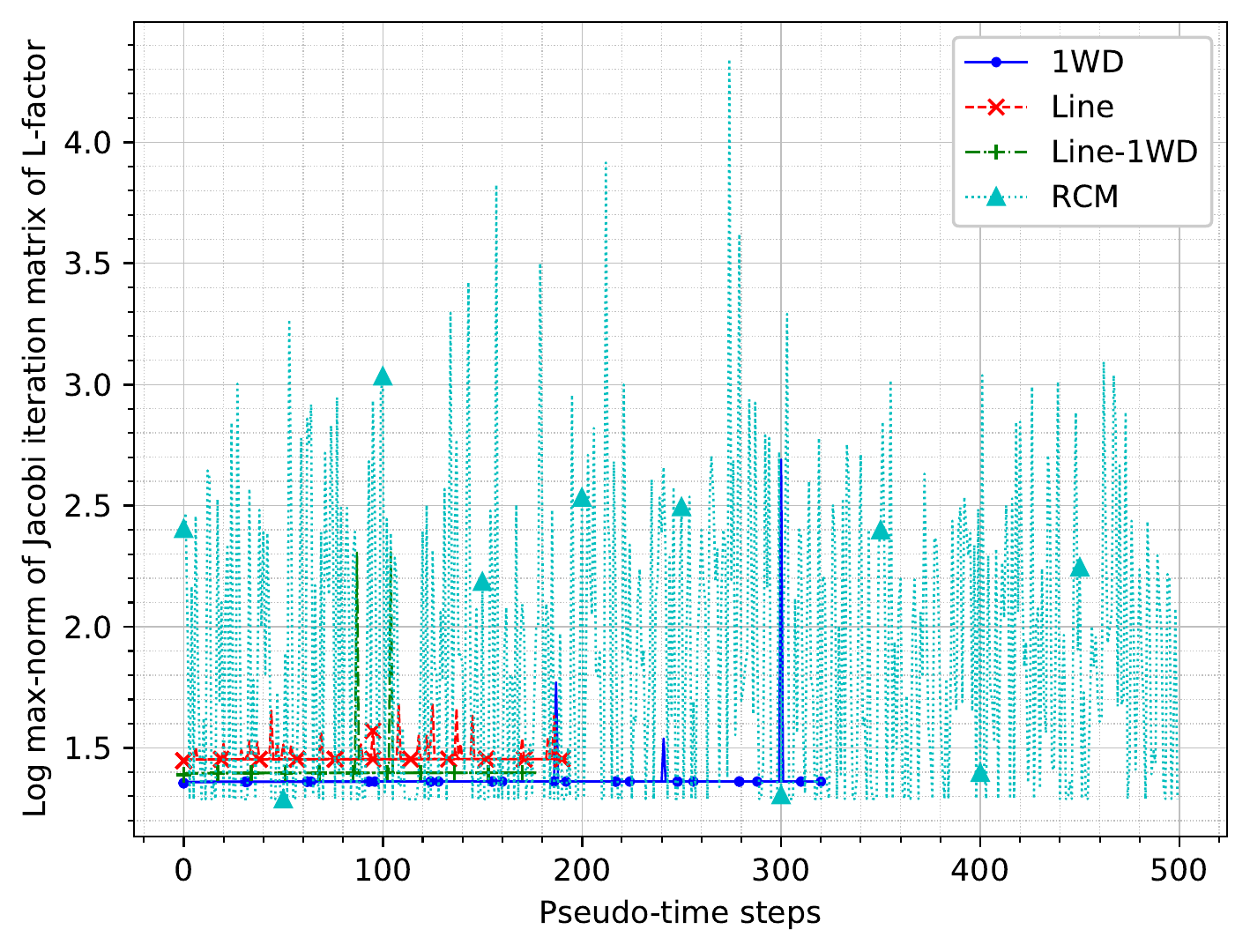}
	\subcaption{Norm of Jacobi iteration matrix for $L$-factor}
    \label{fig:ldom-abilu-noscale-orderings}
	\end{subfigure}
	\begin{subfigure}{0.52\linewidth}
    	\centering
	    \includegraphics[scale=0.52]{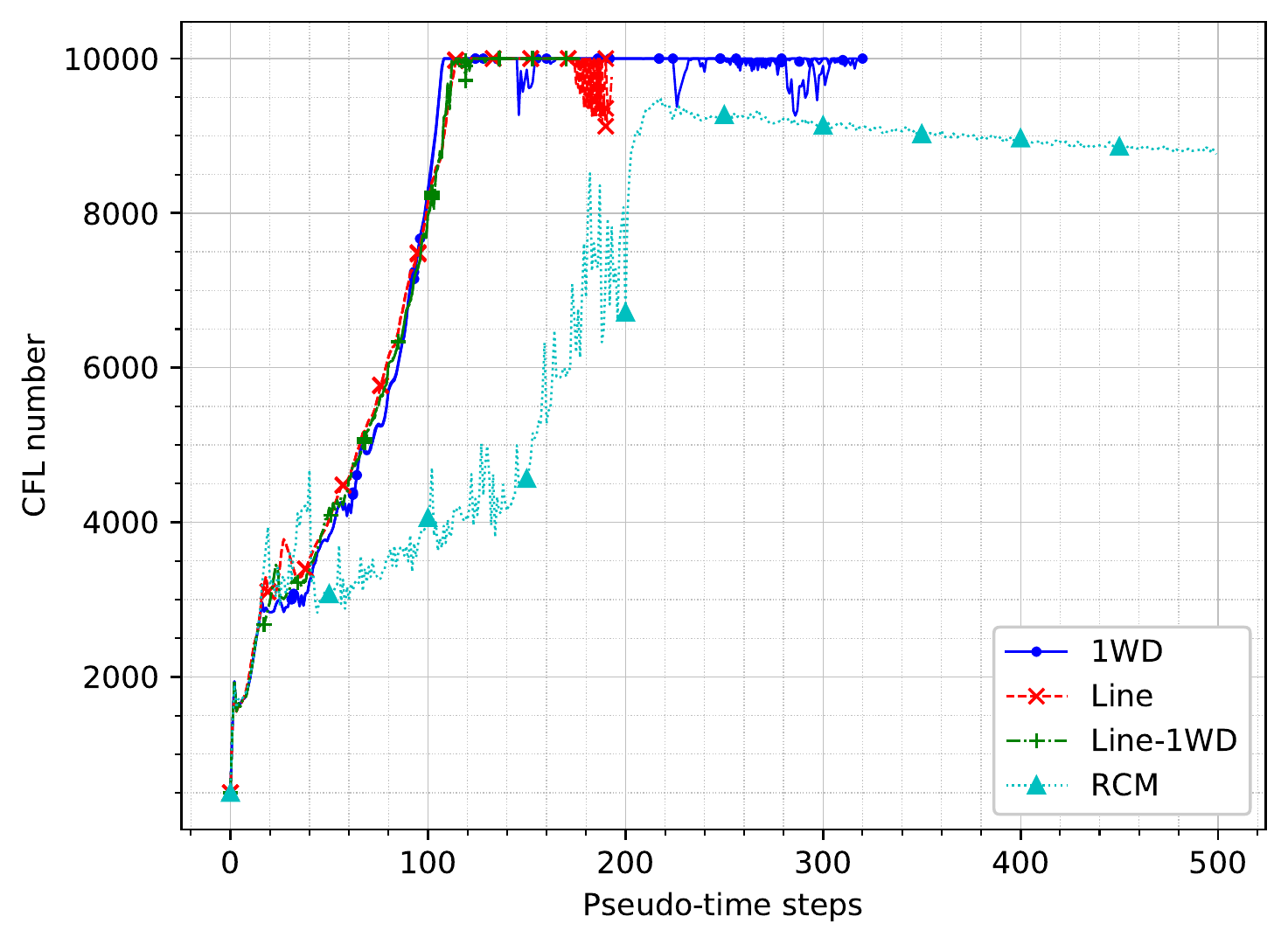}
    	\subcaption{CFL number}
        \label{fig:cfl-abilu-noscale-orderings}
	\end{subfigure}
    \caption{Properties of async. BILU factorization, and CFL number, w.r.t. pseudo-time steps (62 threads)}
    \label{fig:ordering-precinfo}
\end{figure}

Again, we can look to the properties of the $L$ and $U$ factors to attempt to understand why certain orderings work better. Figure \ref{fig:ordering-precinfo} shows properties of the asynchronous block ILU factorization for the case of 62 threads.
The RCM ordering leads to a much less accurate (in the max vector norm) ILU factorization than any of the other orderings, while the hybrid line-1WD ordering results in the most accurate ILU factorization (figure \ref{fig:ilures-abilu-noscale-orderings}).
Next, even though the RCM ordering leads to a $U$ that has slightly lower iteration matrix norm \emph{on average} (figure \ref{fig:udom-abilu-noscale-orderings}), it has sharp spikes at some time steps. For all the orderings the Jacobi iteration matrix norm of the $U$-factor rises as the CFL number increases (figure \ref{fig:cfl-abilu-noscale-orderings}). Once the CFL number reaches its limit of 10,000, the norm stabilizes and remains approximately constant.
Finally, the norm of the iteration matrix for the $L$ factor is clearly very unfavourable in case of the RCM ordering when compared to others (figure \ref{fig:ldom-abilu-noscale-orderings}), owing to very high peaks reached at many of the time steps.

%

\begin{figure}[!htb]
  \begin{subfigure}{0.5\linewidth}
    \centering
    \includegraphics[scale=0.55]{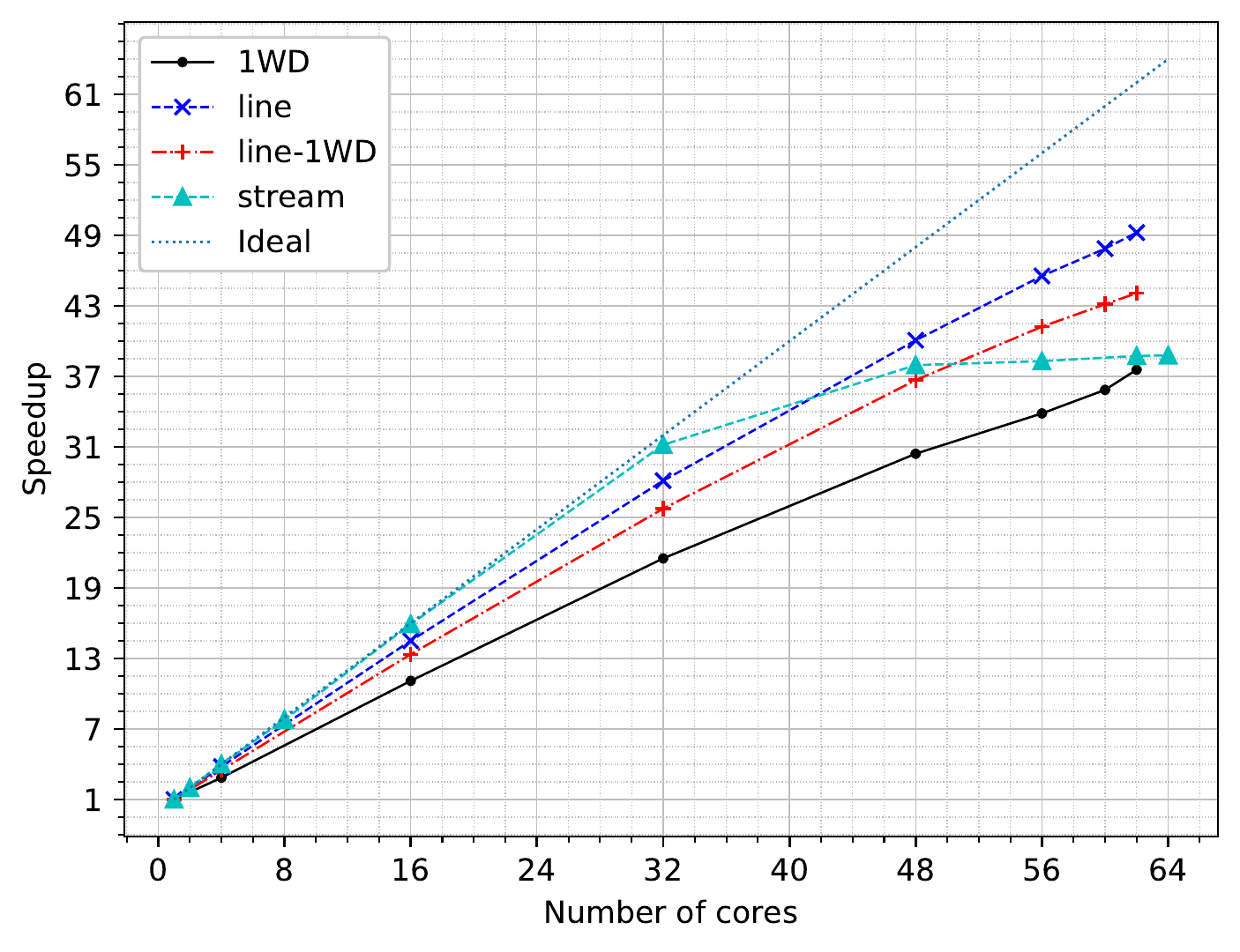}
    \subcaption{Speedup}
    \label{fig:orderings-speedup}
  \end{subfigure}
  \begin{subfigure}{0.5\linewidth}
    \centering
    \includegraphics[scale=0.55]{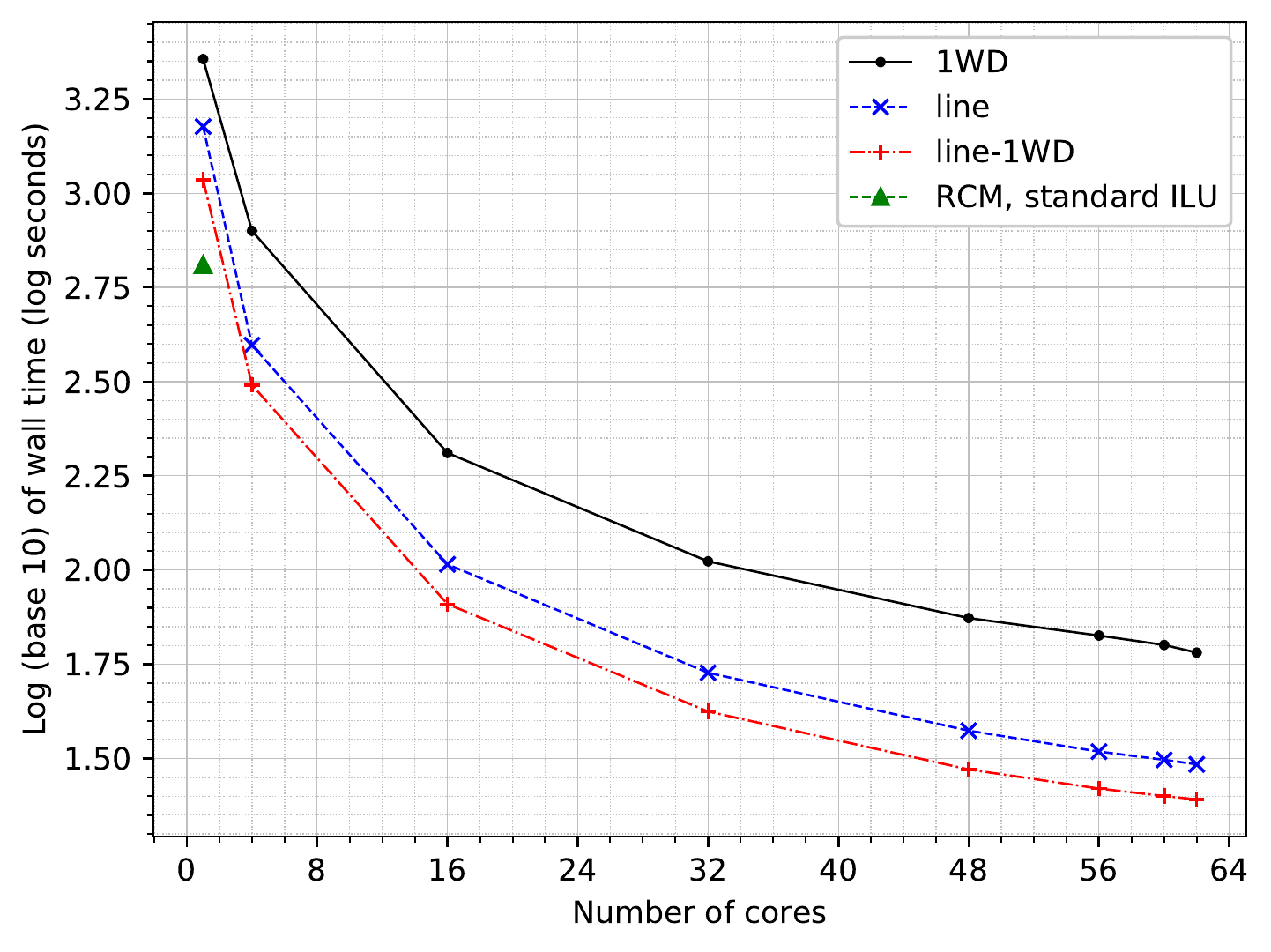}
    \subcaption{Wall clock time}
    \label{fig:orderings-walltime}
  \end{subfigure}
  \caption{Performance of asynchronous block ILU using different orderings (without scaling); though the line ordering scales best, the line-1WD ordering is the fastest up to 62 cores}
  \label{fig:orderings-performance}
\end{figure}

Finally, we show the performance in terms of wall-clock time and strong scaling. As before, we emphasise that data points represent the cumulative time taken by all preconditioning operations over the entire non-linear solve, while other operations are excluded.
The ordering has a large impact on speedups (shown in figure \ref{fig:orderings-performance}), as expected from convergence in terms of FGMRES iterations seen in figure \ref{fig:abilu-orderings-resconv}. Though the line ordering scales best (figure \ref{fig:orderings-speedup}), the line-1WD ordering is actually the fastest (figure \ref{fig:orderings-walltime}). We include one data point for the standard sequential block-ILU preconditioner with RCM ordering in figure \ref{fig:orderings-walltime} for context.

In the speedup plot (figure \ref{fig:orderings-speedup}), we have included the scaling of the Stream benchmark. The scaling of the line-based orderings is comparable to Stream scaling up to a moderate number of cores. As in the inviscid case, the better asynchronous ILU variants continue to scale to higher core counts than Stream, which indicates that they are not quite as highly bandwidth-limited.
 
In the interest of reproducibility and transparency, our codes are available online under the terms of the GNU General Public License. The asynchronous iterations are implemented as a separate library \footnote{\url{https://github.com/Slaedr/BLASTed}, commit 8522c6d26d21c8e63cc09761dd5b1362258a7b6a}. The finite volume CFD solver is available as a code \footnote{\url{https://github.com/Slaedr/FVENS}, commit d05df7fd2b85fd61434dce973e275ccc9d6ea80e} that optionally links to the library. 


\section{Conclusions}
Convergence proofs for asynchronous ILU iteration have been extended to the case of asynchronous block ILU iteration; the latter shows a better theoretical convergence property. From the numerical results, it can be seen that asynchronous block ILU factorization and asynchronous forward and backward block triangular iterations hold promise for fine-grain parallel solution of compressible flow problems. For this coupled system of PDEs, the block variants are much more effective than the original asynchronous ILU preconditioning and scalar asynchronous relaxation for triangular solves.

A second conclusion that can be drawn from this work is that typical grid orderings used for sequential ILU preconditioners may not yield satisfactory results for asynchronous ILU preconditioners, especially for viscous flows. A hybrid `line-X' ordering scheme has been introduced to deal with this issue. For external aerodynamics with viscous flow, one-way dissection (1WD), line and hybrid line-1WD orderings are seen to be good candidates. Further studies are needed on why the 1WD ordering works better than the reverse Cuthill-McKee (RCM) ordering for asynchronous ILU for these problems, even though it is inferior for sequential ILU factorization. In addition, further efforts are needed to develop an ordering and memory storage layout that work well on graphics processing units.

It remains to demonstrate application to larger and more complex problems of compressible flow. Further, these asynchronous iterations would be very useful as fine-grain parallel multigrid smoothers. In this regard, an investigation of the smoothing property of asynchronous block ILU iterations would be useful. Another line of work would be to extend the ideas presented here to parallel threshold ILU factorization \cite{ilu:parilut}, which would be useful for solving some types of linear problems to deeper convergence.

\printbibliography

\end{document}